\documentclass[12pt]{article}

\usepackage{graphicx}
\usepackage{amsmath,amssymb}
\usepackage{graphics}
\usepackage{graphicx}
\usepackage[latin1]{inputenc}
\usepackage{epsfig}
\usepackage{float}

\newcommand{\R}{\mathbb{R}}

\begin{document}
%
\title{{High-Order Discretization of Backward\\ Anisotropic Diffusion\\ and Application to Image Processing}}
\author{\kern-0.0truecm
\normalsize{ Lorella Fatone} \\
{\kern-0.0truecm\small Dipartimento di Matematica, Universit\`a di Camerino,}\\
{\kern-0.0truecm\small  Via Madonna delle Carceri 9, 62032
Camerino (MC),
Italy, lorella.fatone@unicam.it}\\[4mm]
\normalsize {\kern-0.0truecm Daniele Funaro}\\
{\kern-0.0truecm\small Dipartimento di Scienze Chimiche e Geologiche, Universit\`a  di Modena e Reggio Emilia, }\\
{\kern-0.0truecm\small   Via Campi 103, 60121 Modena (MO), Italy, daniele.funaro@unimore.it}\\
}
%
%
%
%
%
%
%
%
%
\date{}
\maketitle
\begin{abstract}
Anisotropic diffusion is a well recognized tool in digital image processing,
including edge detection and denoising. We present here a particular nonlinear time-dependent operator 
together with an appropriate high-order discretization for the space variable.
The iterative procedure emphasizes the contour lines encircling the objects, paving the way to accurate reconstructions 
at a very low cost. One of the main features of such an approach is the possibility of relying on a rather
large set of invariant
discontinuous images, whose edges can be determined without introducing any approximation.
\end{abstract}


\section{Introduction}\label{sec1}

The literature on digital image processing is rather rich. For this reason, we restrict our reference list by only mentioning a few publications: \cite{Gonzalez}, \cite{Granlund},  \cite{Sapiro}, \cite{Solomon}, \cite{WeickertBook}.
For at least three decades, partial differential equations (PDEs) have been contributing to the progress in image processing. 
Since the seminal works of Mumford-Shah \cite{Mumford}  and Perona-Malik  \cite{Perona}, variational and PDE-based methods have found a wide range of applications.
Some works
in these directions, particularly related to anisotropic diffusion using PDEs, are  for instance:  \cite{Alvarez}, \cite{Lions}, \cite{Buades},   \cite{Felsberg}, \cite{Lindenbaum},  \cite{Weickert}.
Here, images are analyzed in a continuous context through appropriate nonlinear time-dependent equations,
suitably discretized. The equations are dissipative in nature, and the nonlinear coefficients
act with different strengths, depending on the intensity of the local gradient. 
In this framework, for example, edge detection algorithms tend to subdivide,
during time evolution, the image into uniform regions, with the aim of providing successively
a reasonable recovery of the contour lines of the objects involved.
We use a similar approach in this paper, realizing some meaningful improvements
with respect to the more classical versions. First of all, the discretization of the PDE
will be of a higher-order. This property is achieved by using a larger stencil and, in addition to the standard
grid-points, an intermediate grid (the so-called ``staggered grid") built on half-fractional indices.
Such a scheme can be adopted in several other situations, not necessarily related to the specific
cases studied here.

The second change we propose concerns the nonlinear diffusion coefficient, which is 
built in order to ensure a large space of ``invariant images". These are images
presenting sharp discontinuities that are not modified by the action of the operator. The algorithm
of edge detection in these situations does not introduce any approximation, so that the
contours remain perfectly represented. Unfortunately, the differential equation resulting
from this modification is rather dissipative in presence of sharp variations
(with the exception of the above mentioned invariant images). This is in contrast
with the usual approach and makes the algorithm non competitive when
iterating many times with a time-step chosen according to the standard 
arrow of time.
Things are instead getting better if we move backwards (negative time).
Reverse nonlinear diffusion is in general an unstable process. This is however true
for the continuous problem. In truth, we will implement the discrete scheme for just one or a
few steps, so that the method is not going to be a source of instabilities.
In this way, with a proper triggering of the parameters, the variations tend
to become sharper,  helping the process of  edge determination.
Some type of noise (such as  ``salt and pepper") is also greatly amplified, 
so that it can be easily detected and successively 
filtered in an appropriate fashion. For the sake of brevity, this last aspect will  not be considered
in the present paper.

The applications presented here are very preliminary. Improvements can be
certainly examined, in conjunction with a series of other specialized filtering
procedures. Actually, a  single step of the method is very low-cost and can provide interesting
information about the image under study. For example, this prediction can be used as a
starting guess in view of implementing other, more sophisticated, techniques, 
as those already available in the market.

\section{Introducing the algorithm in 1D}\label{sec2}

Let $N$ denote the number of pixels along a given interval.
Most of the times, one has $N=256$.
Their distance will be denoted by $h$.  In practice, given $x_0$, we will deal
with the points:
\begin{eqnarray}\label{nodiD}
x_{i+1}=x_{i}  + h, \qquad\qquad  i=0,1,\ldots N,
\end{eqnarray}
where we added two extra pixels at the beginning and at the end.
We will also consider the $N+3$ shifted points, obtained from the expression:
\begin{eqnarray}\label{nodiDshif}
x_{-1/2}=x_{0}-\dfrac{h}{2}, \qquad\qquad 
x_{i+1/2}=x_{i}+\dfrac{h}{2}
, \qquad i=0,1,\ldots N+1.
\end{eqnarray}
We start from our image, which is represented by the values $u_{i}, i=1,\ldots N$, belonging to
a palette with tones of grey from 1 to $K$. In typical applications
we have $K=256$. The minimum grey level is 1 (black) and  the maximum  level   $K=256$ (white). 
We extend the image outside its limits by defining
$u_0=u_1$ and $u_{N+1}=u_N$. This will correspond later to what we call Neumann boundary conditions.

%
%
We first build the (first-order) variations:
\begin{equation}\label{diff}
d_{i+3/2}=\frac{1}{h}\left(u_{i+2}-u_{i+1}\right)\approx \frac{du}{dx}(x_{i+3/2}),
\end{equation}
where $u$ is some differentiable function assuming the values $u_i$ at the points $x_i$.
Successively, we introduce an approximation of the second
derivative operator:
\begin{eqnarray}\label{secondD}
s_{i+3/2}&=&\frac{1}{h}\left(
\frac{u_{i+3}-u_{i+1}}{2h}-\frac{u_{i+2}-u_{i}}{2h}\right)= \frac{1}{2h^2}\left(u_{i}-u_{i+1}-u_{i+2}+u_{i+3}\right)\nonumber\\
&=&
 \frac{1}{2h}\left(d_{i+5/2}-d_{i+1/2}\right) \approx
\frac{d^2u}{dx^2}(x_{i+3/2}), \qquad i=0,\ldots, N-2.
\end{eqnarray}
It is possible to show that the above formula is exact when $u$ is a polynomial
of degree less than or equal to three.

We then introduce the quantities:
\begin{equation}\label{Rquant}
\displaystyle R_{i+3/2}=  \sqrt{\frac{(d_{i+3/2})^2}   {1+(d_{i+3/2})^2}} \ s_{i+3/2},
, \qquad  i=0,\ldots, N-2, 
\end{equation}
where we also define:
\begin{equation}\label{neum}
 R_{-1/2}=R_{5/2}, \ \    R_{1/2}=R_{3/2}, \ \    R_{N+1/2}=R_{N-1/2}, \ \    R_{N+3/2}=R_{N-3/2}. 
\end{equation}

The next step is to interpolate the above values from the fractional to the integer grid. To this purpose we adopt the formula:
\begin{equation}\label{intr}
 R_{i}\approx -\frac{1}{16} R_{i-3/2}+\frac{9}{16} R_{i-1/2}+\frac{9}{16} R_{i+1/2}-\frac{1}{16} R_{i+3/2},
\end{equation}
based on a Lagrange type interpolation. This is also exact up to polynomials of degree
three.

The final step is to update the vector $ u_{i}$, $i=1,\ldots N$. This is done as follows:
\begin{equation}\label{step}
 u_{i} \ \leftarrow \  u_{i} + \gamma  R_{i}, \quad\quad i=1,\ldots N,
\end{equation}
where $\gamma \in {\bf R}$ is a suitable parameter to be determined later.
Again, we can prolong the set of values by defining
$u_0=u_1$ and $u_{N+1}=u_N$. In this way, if necessary, we can restart
the cycle from formula (\ref{diff}) and get a successive update of the values
$u_i$, $i=1,\ldots N$. In the next section we explore the reasons of the above construction.

\section{Preliminary considerations in 1D}\label{sec3}

There are strong connections with partial differential equations.
Indeed, let $u=u(t,x)$ be the solution of the following  time-dependent nonlinear problem
defined in a finite interval:
\begin{equation}\label{equa}
\frac{\partial u}{\partial t}= \sqrt{\frac{d^2} {1+d^2}}\ \frac{\partial^2 u}{\partial x^2},
\end{equation}
where $d=\partial u/\partial x$. The problem is supplied with initial conditions
and boundary constraints of Neumann type at the extremes of the interval.

The construction proposed in the previous section corresponds to a
finite-differences approximation of the above  time-dependent nonlinear problem, with respect to the space variable.
This is followed by a first-order time discretization using the explicit Euler
scheme. One step of the algorithm is actually represented by the
passage in (\ref{step}). Note that the parameter $\gamma$ in (\ref{step}) may assume
positive or negative values. This means that time can flow in both directions,
so that our differential problem in not necessarily dissipative in nature.
We are not interested to study the behavior of (\ref{equa}) in time,
and, in the applications that follow, we will implement (\ref{step}) just for
one or a few  steps.

Contrarily to classical methods (see, for example,  section 5), the anisotropic coefficient $\sqrt{d^2/(1+d^2)}$ 
in  (\ref{equa}) is less pronounced when the derivatives are small, therefore the solution
tends to smooth out the edges during the time  evolution. For this reason, we will prefer to work
with negative values of $\gamma$. Running backwards brings notoriously to ill-posed problem,
although the question can be handled with the due hypotheses (see, e.g., \cite{Hummel}).
We will use, however, very few steps of the algorithm and the parameters will be calibrated
in order to stay away from impracticable situations.



We would like to mention some features of the algorithm in order to justify
the choices we made.
Suppose that between the given points $x_i$ and $x_{i+1}$ of a given 1D image we have a  jump.
With this we mean that $u_j=V_1$ for $j\leq i$, and $u_j=V_2$ for $j\geq i+1$, where 
$V_1, V_2 \in  \R$, with $V_1\not =V_2$. Then, after a cycle of  (\ref{step}), these values remain unchanged.
This is true because the variations in (\ref{diff}) are all zero with the exception of $d_{i+1/2}$. On the other hand, it easy to check that
$s_{i+1/2}=0$ (see (\ref{secondD})). Since, for any $k$, the values $R_{k+1/2}$ in  (\ref{Rquant}), (\ref{neum}) contain the products
$d_{k+1/2}s_{k+1/2}$, they turn out to be identically zero. Thus, the
algorithm exactly preserves isolated discontinuities, and therefore
does not introduce any dissipation, whatever is $\gamma$. It is possible to check that
this property holds true also for distributions where $u_j=V_1$ for
odd indices and $u_j=V_2$ for even indices. There are, in fact, plenty
of basic configurations that are preserved by the algorithm.

\begin{figure}[ht!] 
\vspace{-.4cm}
  \hspace{-.6cm}\begin{minipage}[b]{0.53\linewidth}
    \includegraphics[width=1.1\linewidth ]{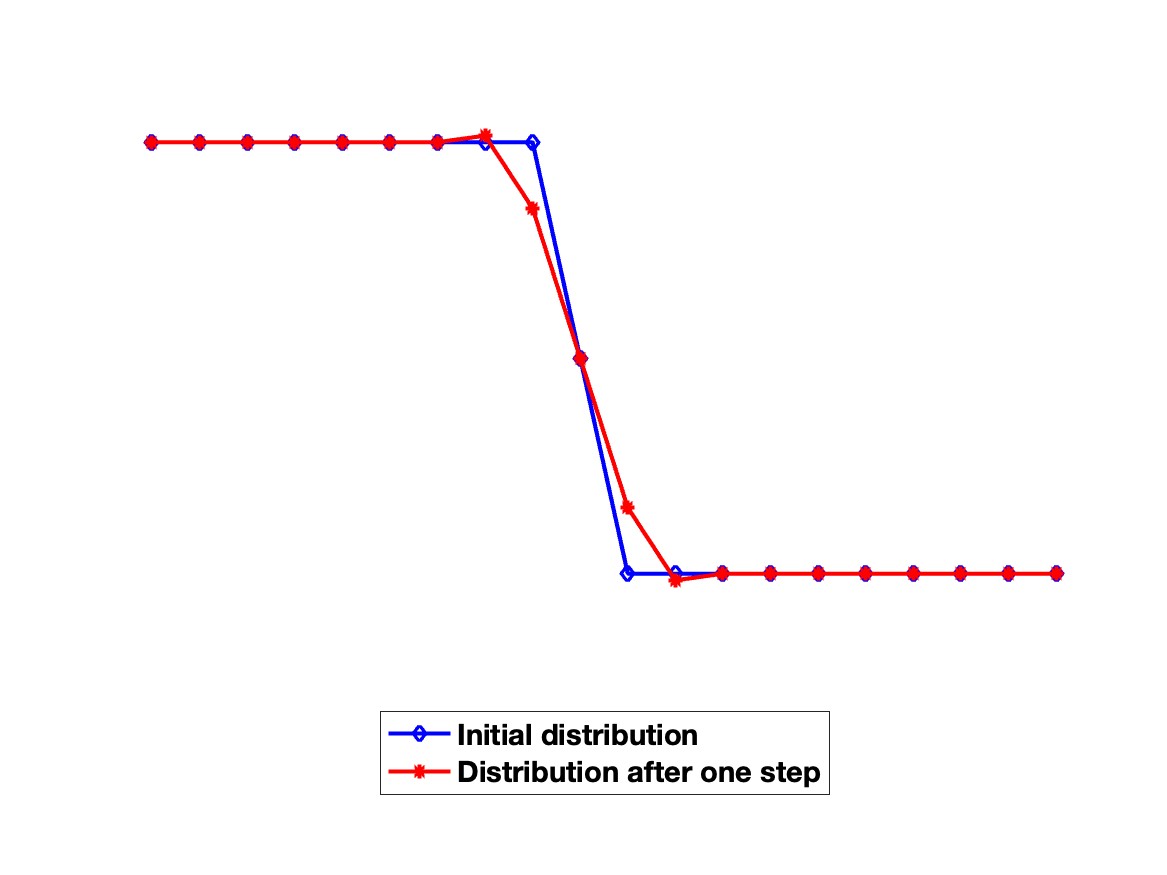} 
  \end{minipage} 
  \hspace{-.8cm}\begin{minipage}[b]{0.53\linewidth}
    \includegraphics[width=1.1\linewidth]{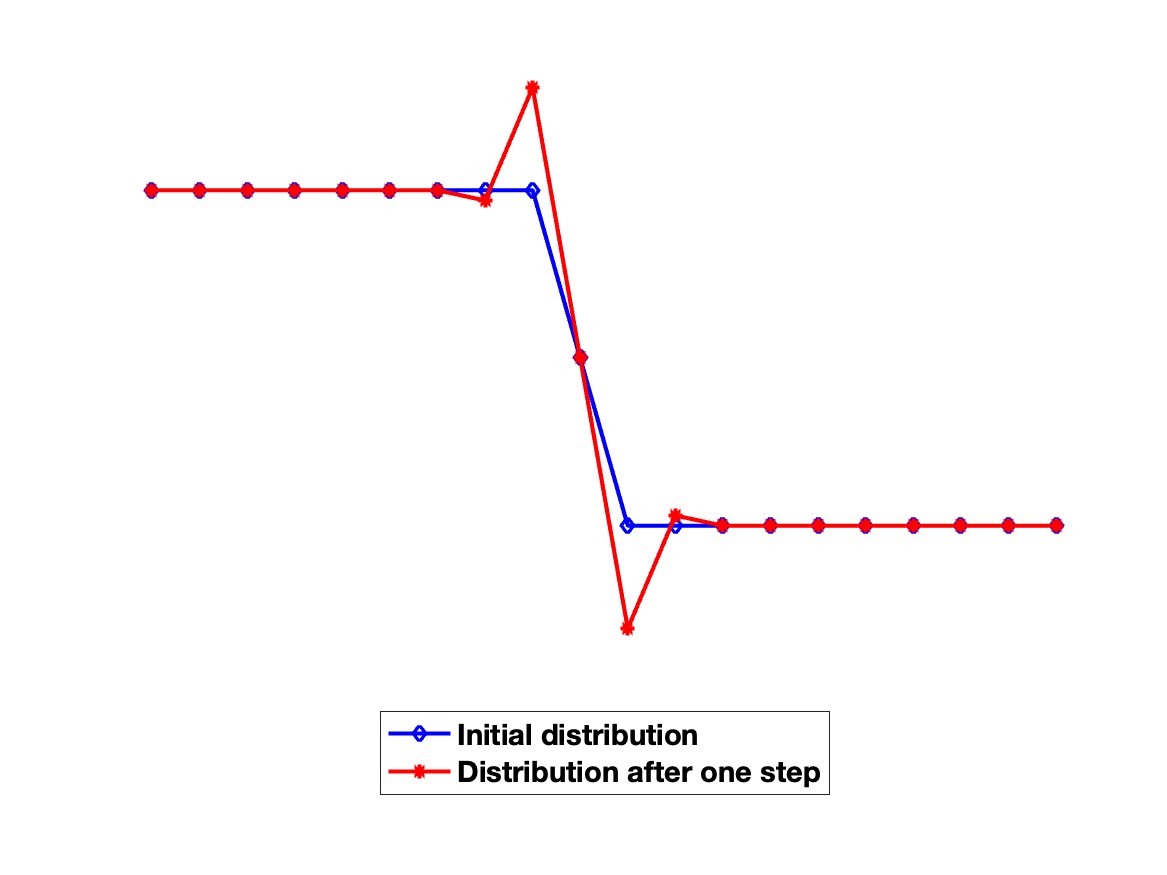} 
  \end{minipage} 
 \vspace{-.3cm}
     \caption{ \small  Transformation of the distribution $u$ corresponding to (\ref{solgradino}) after one step of the algorithm  (\ref{step}) (left: $\gamma>0$, right: $\gamma<0$).
    } 
 \label{fig1} 
\end{figure}

The next interesting case to be studied is the one where we consider the following 1D image:
\begin{equation}\label{solgradino}
u_j=\left\{ \begin{array}{l}
V_1, \quad j \le i-1,\\[3mm] 
\frac12 (V_1+V_2), \quad  j=i,\\[3mm] 
V_2,  \quad j \ge i+1,
\end{array} \right. 
\end{equation}
where  $V_1, V_2 \in  \R$, with $V_1\not =V_2$. In Fig.  \ref{fig1} we show a couple of cases
corresponding to $V_1>V_2$ in (\ref{solgradino}).
When $\gamma$ in (\ref{step}) is positive (plot on the left), after one
step some dissipation is introduced, so that the corresponding picture is smoothed.
When instead $\gamma$ is negative (plot on the right) the jump is emphasized. 
Of course, the parameter $\gamma$ has to be chosen in an appropriate
range (also depending on the values $V_1$ and $V_2$). If $\gamma$ is too small,
nothing interesting happens. Otherwise, if it is too large, the final picture may come out to be heavily
deteriorated. One may also consider to apply a certain number of steps
of this process. The user may then decide the reliability of the outcome
depending on the requested result. Different values of $\gamma$ may be also implemented at each step.

\section{The algorithm in 2D}\label{sec4}
In 2D the natural extension of (\ref{equa}) reads as follows:
\begin{equation}\label{equa2}
\frac{\partial u}{\partial t}= \sqrt{\frac{d^2} {1+d^2}}\ \Delta u,
\end{equation}
where $d=\vert \nabla u\vert$, with  $\nabla u$, $\Delta u$ being the gradient and the Laplace operator,
respectively. Such an equation is defined on a square
and Neumann boundary conditions are imposed. This formulation allows us
to build the algorithm for the two dimensional case, as described here below.

We start by introducing the points:
\begin{equation}\label{nodiD2}
x_{i+1}=x_i + h, \qquad y_{j+1}=y_j+ h,\qquad\qquad  i,j=0,1,\ldots N.
\end{equation}
The values of grey attained at these points are denoted by $u_{i,j}$.
Usually the minimum grey level is 1 and  the maximum grey level is  $256$. 
In this construction we have to assume that the initial picture (only
defined for $i,j=1,\ldots N$) has been bordered in order to satisfy a 
discrete version of Neumann boundary conditions. Note that the values
$u_{0,0}$, $u_{N+1,0}$, $u_{0,N+1}$, $u_{N+1,N+1}$, are irrelevant,
since they will not be used in the algorithm.

\begin{figure}[ht!]%
\centering
\includegraphics[width=0.5\textwidth]{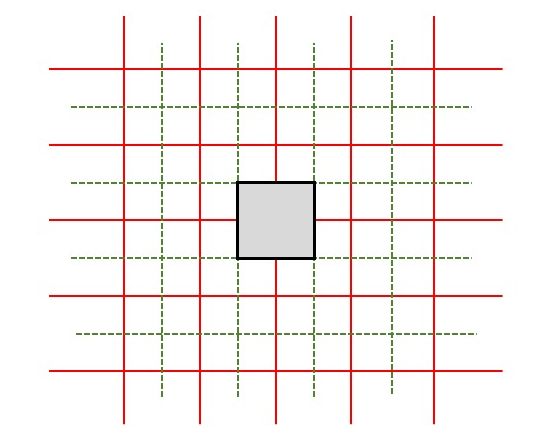}
\vspace{.3cm}
\caption{\small A generic pixel, denoted by the square, is centered at a grid point
corresponding to integer indices (at the crossing of the solid lines). Another grid, with fractional indices, is used
to carry out the computations.}\label{griglia}
\end{figure}

Intermediate grids, where the indices take semi-integer values,
are naturally associated to the primary grid as in (\ref{nodiDshif}), that is:
\begin{eqnarray}\label{nodiDshif2D}
&&
x_{-1/2}=x_{0}-\dfrac{h}{2}, \ \ \
y_{-1/2}=y_{0}-\dfrac{h}{2}, \ \ \
x_{i+1/2}=x_{i}+\dfrac{h}{2}, \ \ \
y_{j+1/2}=y_{j}+\dfrac{h}{2},\nonumber\\[2mm]
&&\hskip7.5truecm i,j=0,1,\ldots N+1.
\end{eqnarray}
The situation is summarized in Fig. \ref{griglia}.

The next step is to define the values at the shifted nodes 
by averaging the values at the integer nodes as follows:
\begin{eqnarray}
&&
u_{i+3/2,j+1}=\dfrac{1}{2}\left(u_{i+1,j+1}+u_{i+2,j+1}\right) , \label{u_nodiDshif2D1}\\
&&
u_{i+1,j+3/2}=\dfrac{1}{2}\left(u_{i+1,j+1}+u_{i+1,j+2}\right) , \label{u_nodiDshif2D2}
\end{eqnarray}
where, for the sake of simplicity, we omitted to specify the range of variability of the indices.
Successively we build the variations in both directions:
\begin{eqnarray}
&&
d^x_{i+3/2,j+3/2}=\dfrac{1}{h}\left(u_{i+2,j+3/2}-u_{i+1,j+3/2}\right) ,\label{var_nodiDshif2D1} \\
&&
d^y_{i+3/2,j+3/2}=\dfrac{1}{h}\left(u_{i+3/2,j+2}-u_{i+3/2,j+1}\right)\label{var_nodiDshif2D2},
\end{eqnarray}
so that the approximation of the Laplace operator ends up to be:
\begin{eqnarray}\label{approxLapl}
l_{i+3/2,j+3/2}&=&\dfrac{1}{2h} \left(  d^x_{i+5/2,j+3/2}-d^x_{i+1/2,j+3/2}+d^y_{i+3/2,j+5/2} -d^y_{i+3/2,j+1/2}\right)\nonumber\\[3mm]
&\approx& (\Delta u)(x_{i+3/2,j+3/2}).
\end{eqnarray}
We then introduce the quantities:
\begin{eqnarray}\label{D_nodiDshif2D}
&&
D_{i+3/2,j+3/2}=\sqrt{\left( d^x_{i+3/2,j+3/2} \right)^2+ \left( d^y_{i+3/2,j+3/2} \right)^2},
\end{eqnarray}
and
\begin{eqnarray}
\displaystyle R_{i+3/2,j+3/2}= \sqrt{\frac{\left(D_{i+3/2,j+3/2}\right)^2}   {1+\left(D_{i+3/2,j+3/2}\right)^2}}\ l_{i+3/2,j+3/2}.\\\nonumber
\end{eqnarray}
We also extend the above values by introducing a double layer in proximity of the boundary,
in order to enforce Neumann type boundary conditions, as done in (\ref{neum}).

As in the one-dimensional case, the successive step is to interpolate the above values from the fractional to the integer grid. We do this by implementing
the Cartesian product of formula (\ref{intr}), so obtaining: 
\begin{eqnarray}\label{fmlaLag}
 R_{i,j}&\approx& \frac{1}{256} 
 \left(  R_{i-3/2,j-3/2}   -9  R_{i-3/2,j-1/2}   -9  R_{i-3/2,j+1/2} + R_{i-3/2,j+3/2}     \right. \nonumber\\
&& -9 R_{i-1/2,j-3/2}   +81 R_{i-1/2,j-1/2}   +81  R_{i-1/2,j+1/2} -9  R_{i-1/2,j+3/2}     \nonumber \\
&& -9 R_{i+1/2,j-3/2}   +81 R_{i+1/2,j-1/2}   +81  R_{i+1/2,j+1/2} -9  R_{i+1/2,j+3/2}     \nonumber\\
&& \left. R_{i+3/2,j-3/2}   -9  R_{i+3/2,j-1/2}   -9  R_{i+3/2,j+1/2} + R_{i+3/2,j+3/2}     \right).
\end{eqnarray}
In order to avoid misunderstandings, let us observe that the number 256 at the denominator
in (\ref{fmlaLag}) represents the algebraic sum of the integer weights
at the numerator, and has nothing to do with the scale of greys of the
color palette.

Finally, an update of the values $ u_{i,j}$ is
realized by setting:
\begin{equation}\label{step2D}
 u_{i,j} \ \leftarrow \  u_{i,j} + {\gamma}  \,  R_{i,j}, \ \quad i,j=1,\ldots N, 
\end{equation}
where again  ${\gamma} \in {\bf R}$ is a suitable parameter that
can be freely chosen depending on the kind of performances  requested.
If another step is required, the procedure restarts from the definitions
in  (\ref{u_nodiDshif2D1}), (\ref{u_nodiDshif2D2}), 
(\ref{var_nodiDshif2D1}), (\ref{var_nodiDshif2D2}). If necessary, different
values of $\gamma$ may be proposed at each step.

\begin{figure}[ht!] 
\hspace{-1.4cm}
  \begin{minipage}[b]{0.5\linewidth}
    \includegraphics[width=1.1\linewidth ]{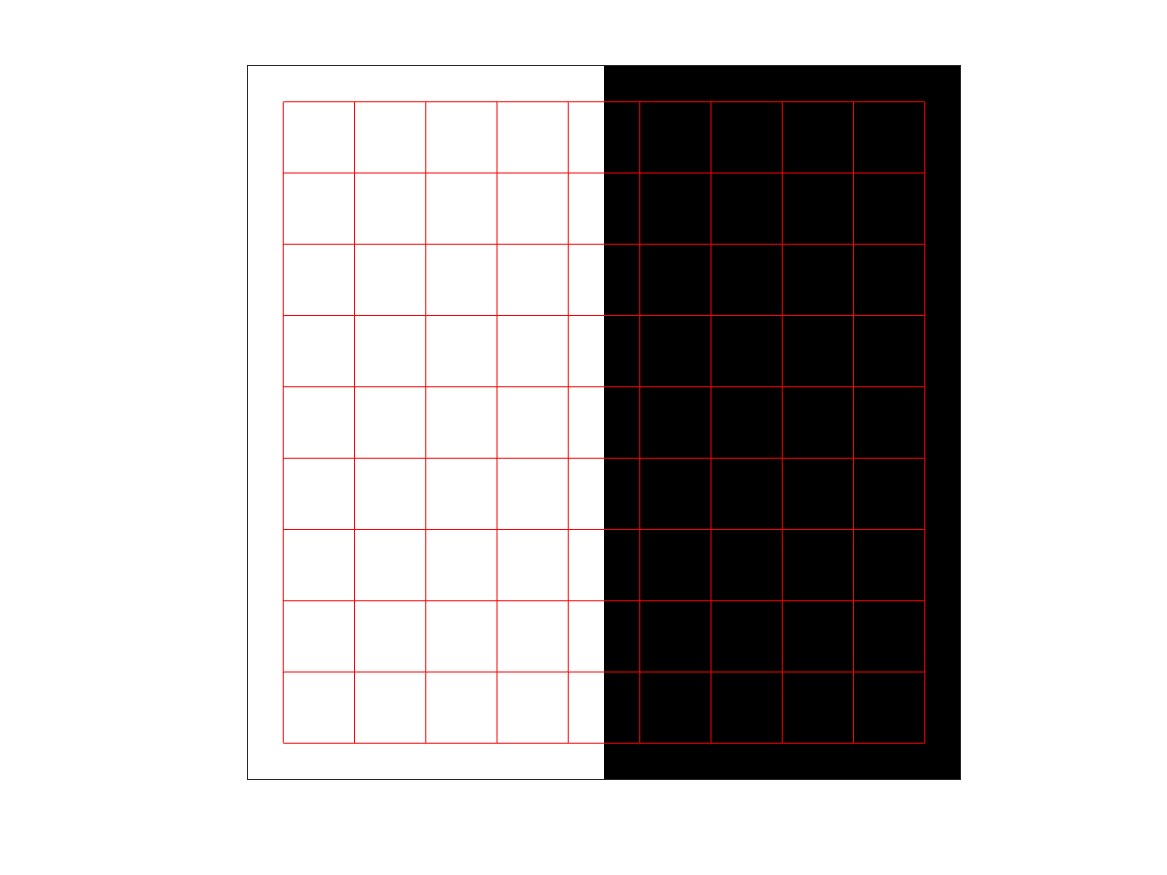} 
  \end{minipage} 
  \begin{minipage}[b]{0.5\linewidth}
    \includegraphics[width=1.1\linewidth]{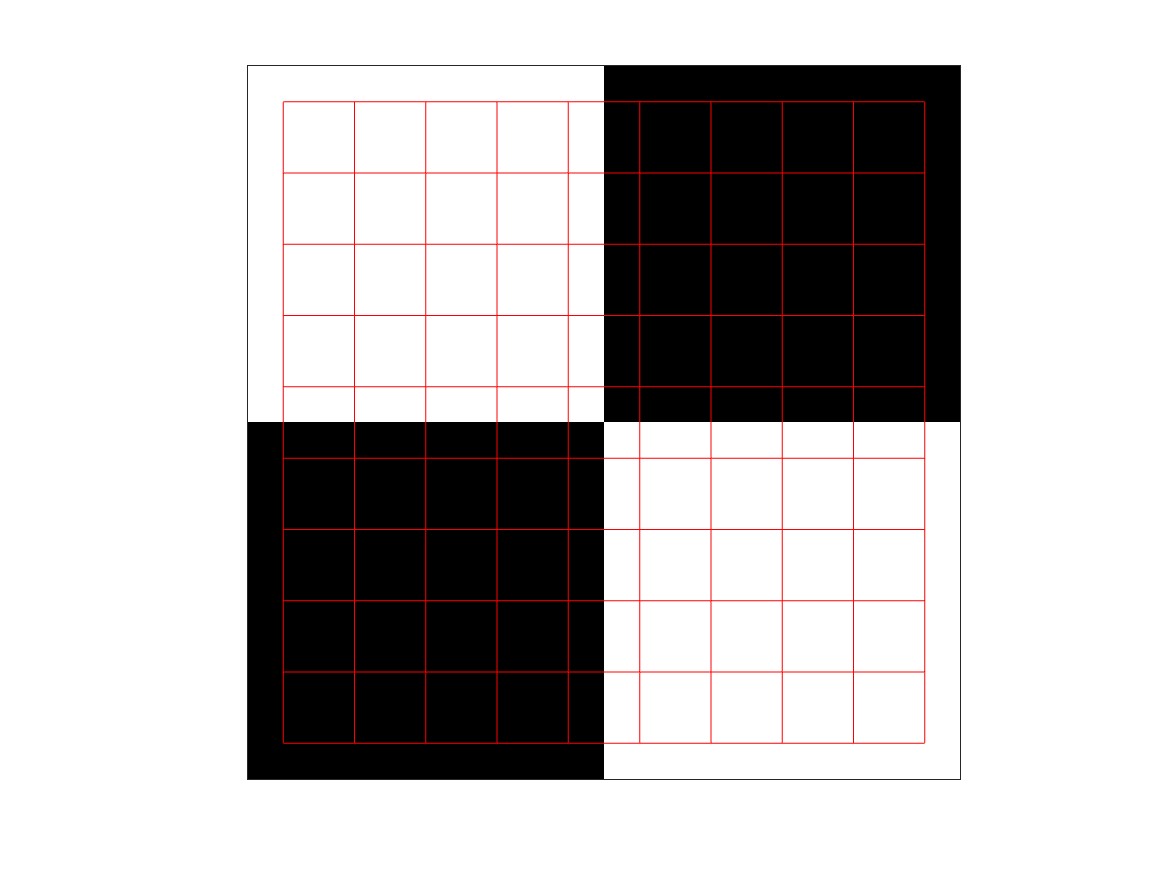} 
  \end{minipage} 
  \begin{minipage}[b]{0.5\linewidth}
  \hspace{-.8cm}
    \includegraphics[width=1.1\linewidth]{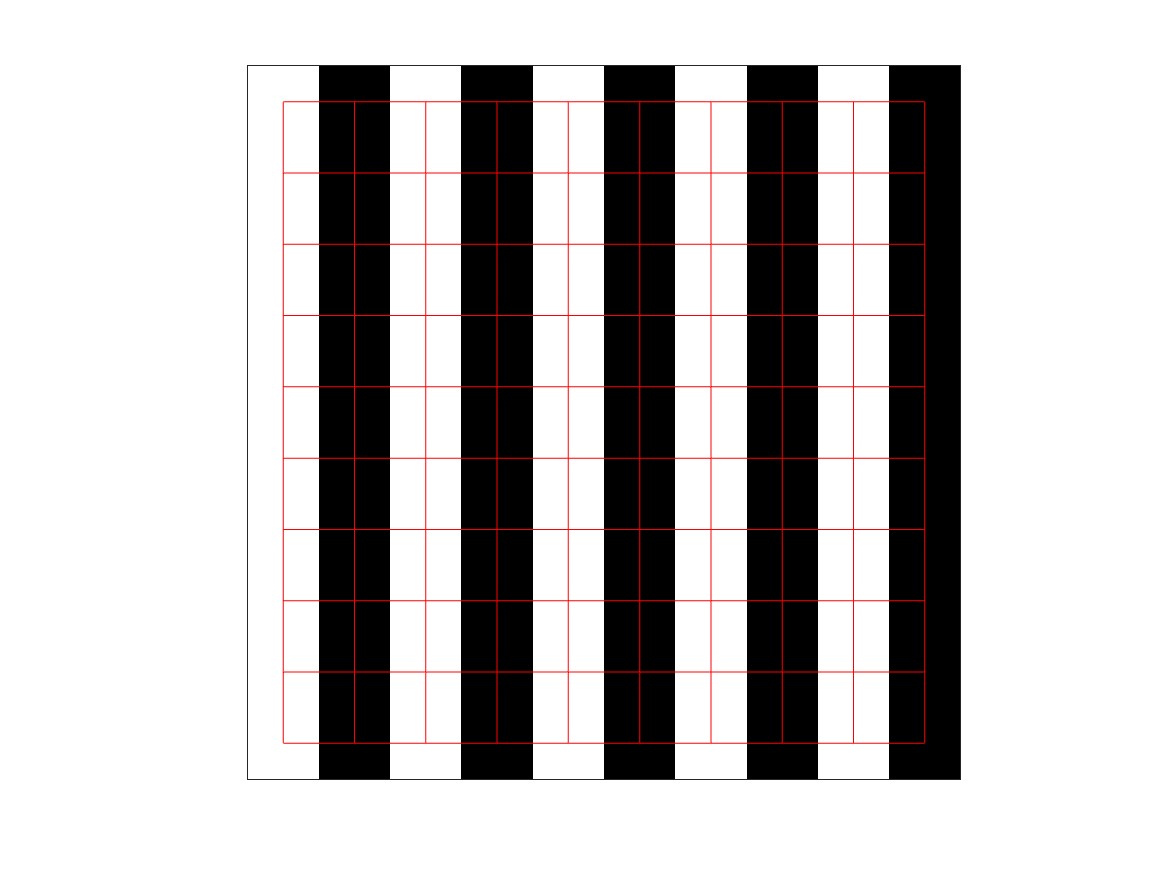} 
  \end{minipage}
  \hfill
  \begin{minipage}[b]{0.5\linewidth}
    \includegraphics[width=1.1\linewidth]{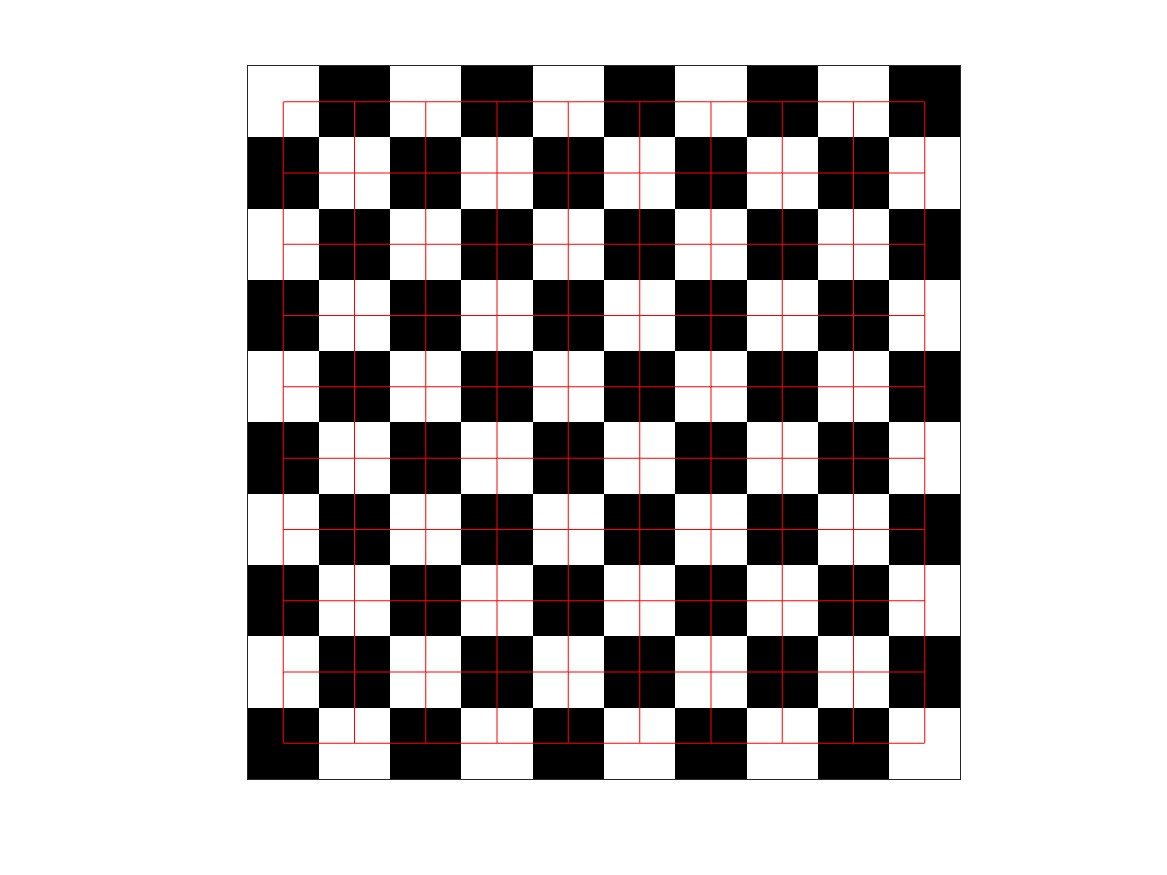} 
  \end{minipage} 
     \caption{\small Although these images present sharp discontinuities, they remain unchanged
     after an application of the nonlinear differential operator resulting from the discretization of the spatial part in (\ref{equa2}).} 
 \label{fig2bw} 
\end{figure}

In Fig. \ref{fig2bw} we display for $N=10$ some black\&white images that
are left unchanged by the algorithm, independently from the value of  $\gamma$ in (\ref{step2D}).
This means that all the values $R_{i,j}$ are zero.
Even if there is the presence of the Laplace operator, the scheme does not
introduce  any diffusion in these cases.

\begin{figure}[ht!] 
\hspace{-1.1cm}
\begin{minipage}[b]{0.37\linewidth}
    \includegraphics[width=1.1\linewidth ]{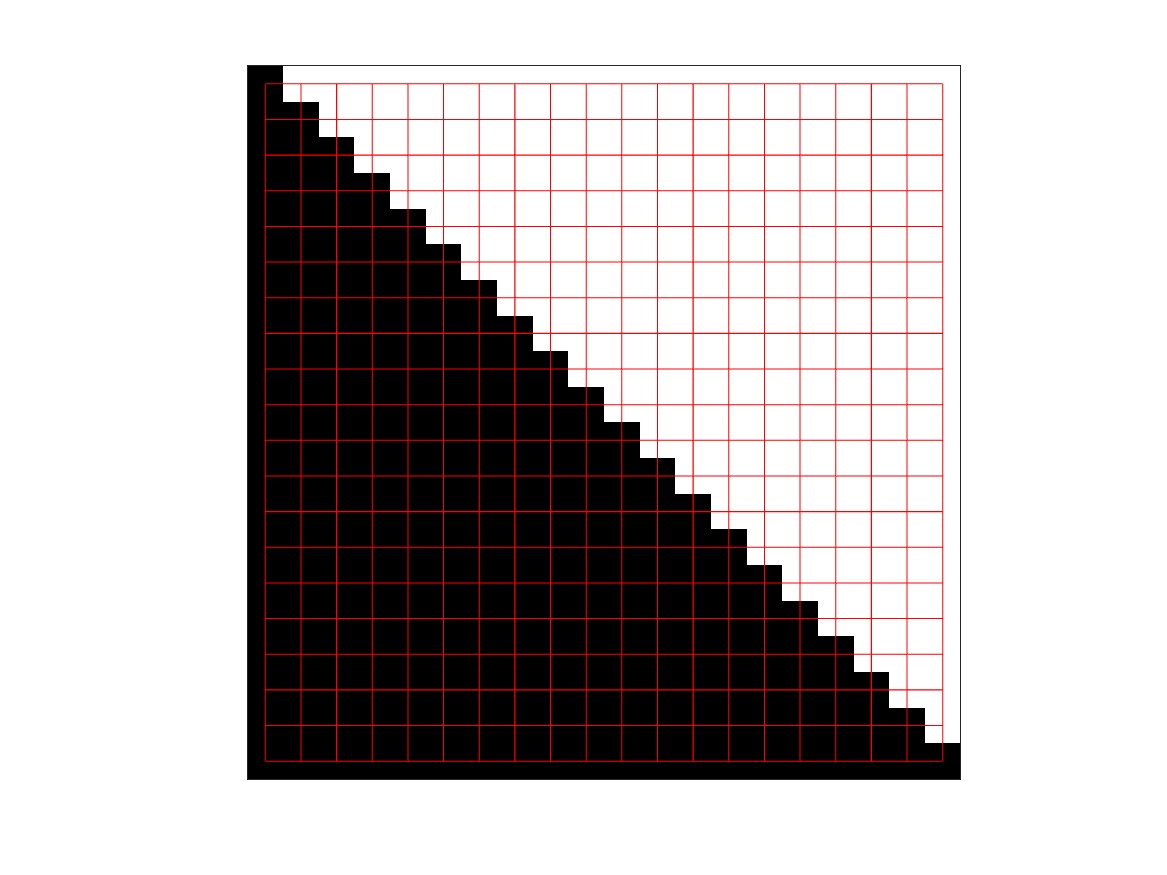} 
  \end{minipage} 
  \hspace{-.5cm}\begin{minipage}[b]{0.37\linewidth}
    \includegraphics[width=1.1\linewidth ]{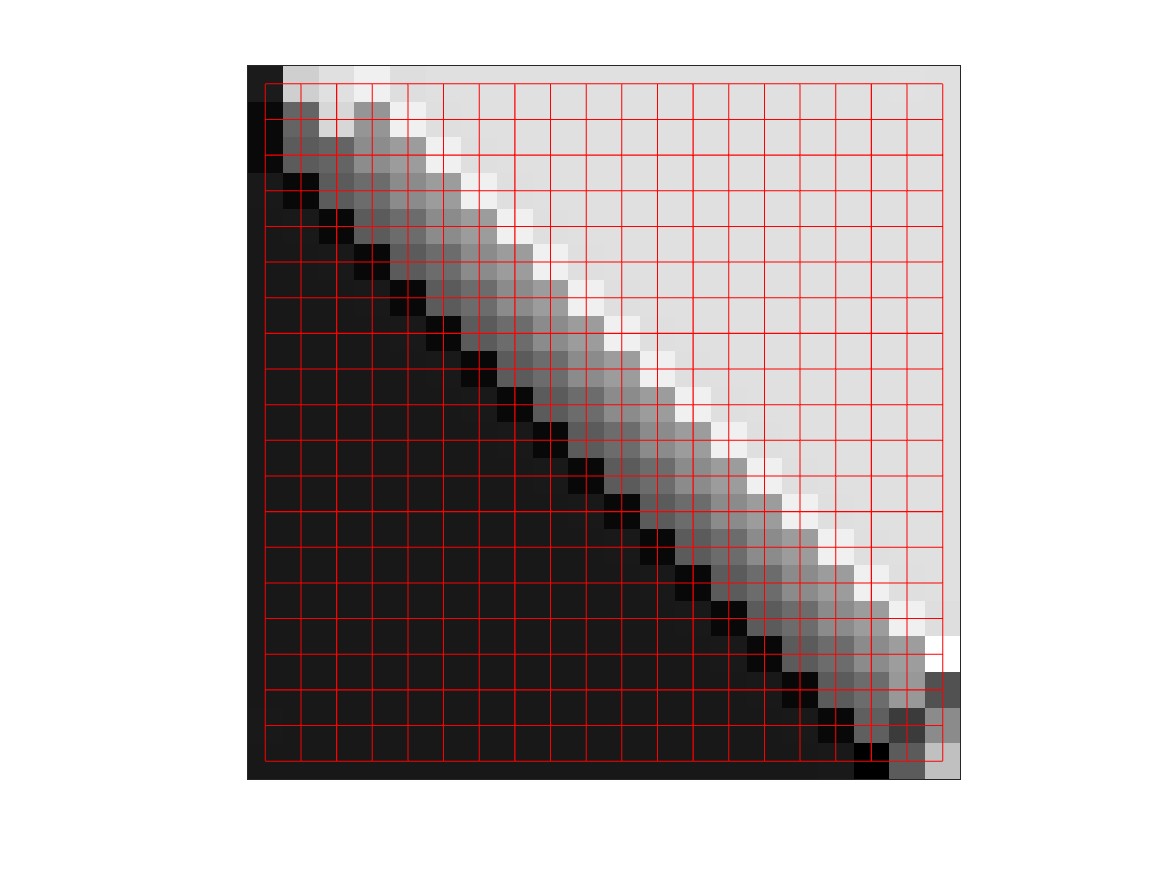} 
  \end{minipage} 
  \hspace{-.5cm}\begin{minipage}[b]{0.37\linewidth}
    \includegraphics[width=1.1\linewidth]{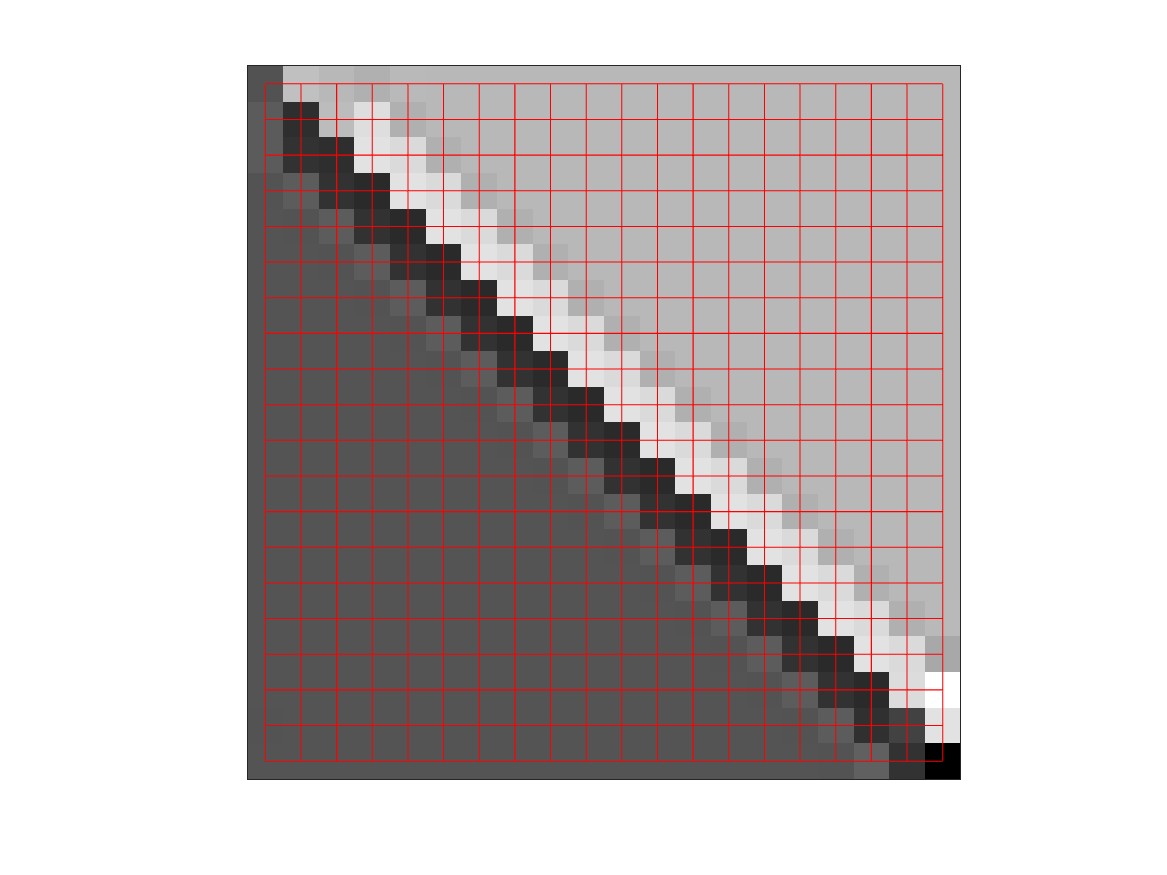} 
  \end{minipage} 
  \vspace{-.3cm}
     \caption{ \small  Transformation of the image on the left  by applying one step of the algorithm (center: $\gamma>0$, right: $\gamma<0$).
    } 
     \vspace{.4cm}
 \label{fig2Diag} 
\end{figure}

Before facing more complex situations (see section 6), a simple preliminary test has been conducted 
on the image of Fig. \ref{fig2Diag}  (left), where we set $N=20$. The pixels corresponding to
$i<j$ are white, whereas the remaining ones are black. This displacement is not
an invariant of the algorithm. Hence, some values $R_{i,j}$ turn
out to be different from zero. After one step of (\ref{step2D}) we are in a situation similar
to that of Fig. 1. For choices of $\gamma$ not excessively large, we get a smoothing of the
picture, if the parameter  $\gamma$  is positive (center image in Fig. \ref{fig2Diag}). 
Indeed, with minor oscillations, the shades change from black to white through a series of intermediate greys.
If, instead, $\gamma$ is negative (image on the right in Fig. \ref{fig2Diag}), we get a kind of overshooting in proximity of the
discontinuity, that amplifies the width of the edge. In these plots, the grey palette has been 
re-scaled in such a way the global minimum is still black and the maximum white.
We also note some discrepancies at the four corners. These are due to the use 
of formula (\ref{fmlaLag}), which requires a certain number of extra values outside the
image domain in order to compute the values $R_{1,1}$,  $R_{1,N}$, $R_{N,1}$, $R_{N,N}$. This says that Neumann type conditions are not always the optimal choice. We do not have, however, other alternatives to suggest at the moment.

\section{A classical scheme for anisotropic diffusion}\label{SecPerona}

The algorithm discussed so far belongs to the category of anisotropic
diffusion schemes. These techniques are based on the  solution
of a nonlinear time-dependent PDE, presenting a dissipative behavior.
The dissipation coefficient generally depends on the magnitude
of the local gradient. 
As an example,  we just introduce here the classical scheme proposed by Perona and Malik in 
their pioneering work \cite{Perona},
which is based on the following nonlinear evolution  equation:
\begin{equation}\label{perona}
\frac{\partial u}{\partial t}= \text{div} \Big(g_{a}(\| \nabla u\|)\ \nabla u\Big),
\end{equation}
where, for a given parameter $a$, the function $g_{a}$ is chosen as follows:
\begin{equation}\label{peronaG}
g_{a} (\eta)=\left(1+\frac{\eta^{2}}{a^{2}}\right)^{\hspace{-.15cm}-1}
\hspace{-.3cm}, \qquad \eta \ge 0.
\end{equation}
%
In this way, the diffusion coefficient is small when  $\| \nabla u\|$  is large (in particular, in proximity of the edges), and  vice versa.  Let us observe that (\ref{equa2}) is not expressed in
divergence form, though this variant could be in principle examined.

We assume that $u(t,x,y)$ in  (\ref{perona}) is defined on a square where homogeneous Neumann type
conditions are enforced at the boundary. The time $t$ belongs to the interval
$[0,T]$, for some $T>0$. The initial configuration is related to the image under study.
%
 After having introduced a time-step $\Delta t>0$, similarly to what has been  done in the previous sections, we denote by $u_{i,j}^{k}$  the pixel representation of our image at the grid points  $(x_{i},y_{j})$ at some given time $t^{k}=k \Delta t$, $k=0,\ldots, K$,  with  $K \Delta t=T$. 
The values are updated according to a classical five-points stencil associated with
the explicit finite difference discretization of
equation  (\ref{perona}), i.e.:  
%
$$
\frac{u_{i,j}^{k+1}-u_{i,j}^{k}}{\Delta t}=\frac1h\left[ G_{i+1/2,j}^{k} \left(  \frac{u_{i+1,j}^{k}  -u_{i,j}^{k}}{h}\right)  - G_{i-1/2,j}^{k} \left(  \frac{u_{i,j}^{k}  -u_{i-1,j}^{k}}{h}\right) + \right. $$
\begin{equation}\label{schemaFD}
\hspace{1.9cm} \left.  G_{i,j+1/2}^{k} \left(  \frac{u_{i,j+1}^{k}  -u_{i,j}^{k}}{h}\right)  - G_{i,j-1/2}^{k} \left(  \frac{u_{i,j}^{k}  -u_{i,j-1}^{k}}{h}\right) \right],
\end{equation}
where: 
\begin{eqnarray}
 G_{i,j}^{k} &=&  g_{a}\left( \sqrt{  \left(  \frac{u_{i+1,j}^{k}  -u_{i-1,j}^{k}}{2h} \right)^2 +    \left( \frac{u_{i,j+1}^{k}  -u_{i,j-1}^{k}}{2h}\right)^2}\right)  \nonumber\\[3mm]
&\approx& g_{a}\left(  \left \|      \nabla  u(t^{k},x_{i},y_{j})\right\|\right), 
  \end{eqnarray}
 and, for instance:
  \begin{equation}\label{diffusionG}
 G_{i+1/2,j}^{k} = \frac12 \left(  G_{i,j}^{k} + G_{i+1,j}^{k}    \right),
  \end{equation}
with similar definitions adopted for the other terms in  (\ref{schemaFD}).

\begin{figure}[ht!] 
\vspace{-.1cm}
\hspace{-.8cm}
  \begin{minipage}[b]{0.5\linewidth}
    \includegraphics[width=1.1\linewidth ]{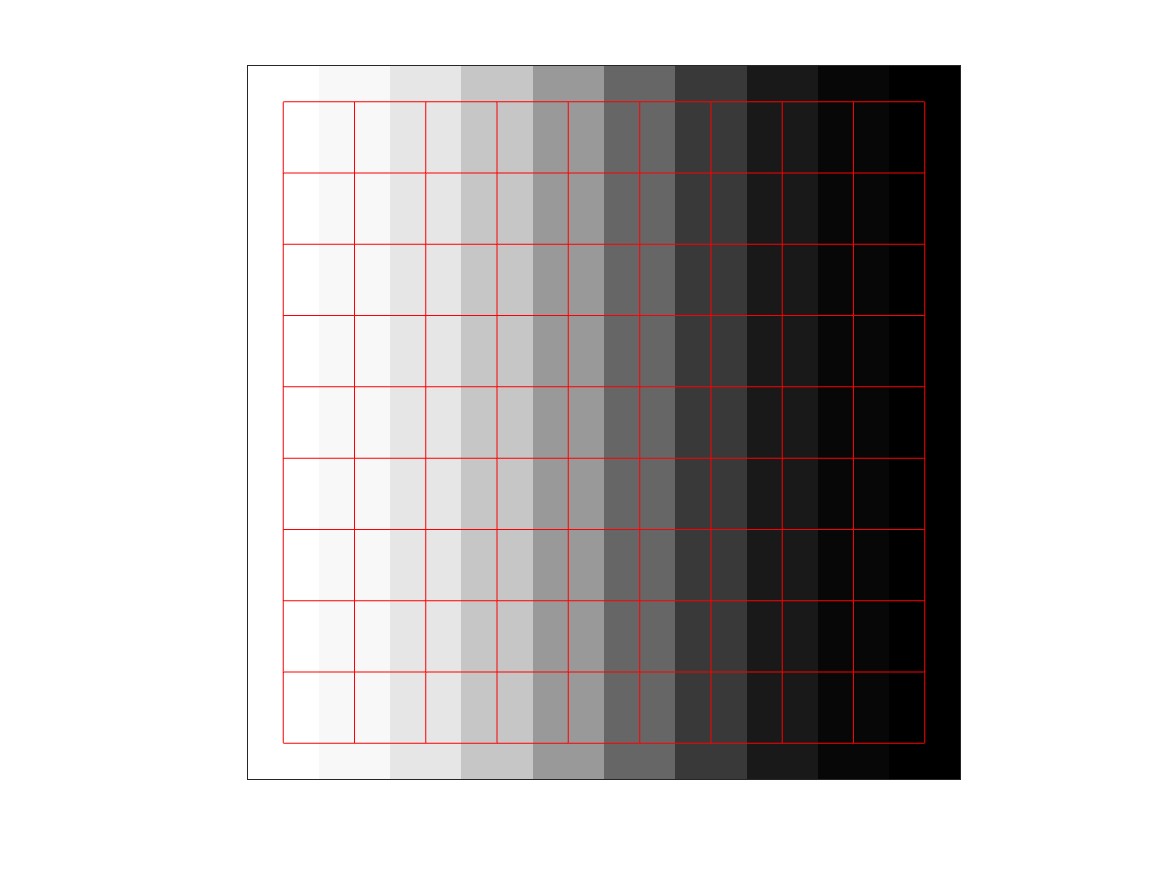} 
  \end{minipage} 
  \begin{minipage}[b]{0.5\linewidth}
    \includegraphics[width=1.1\linewidth]{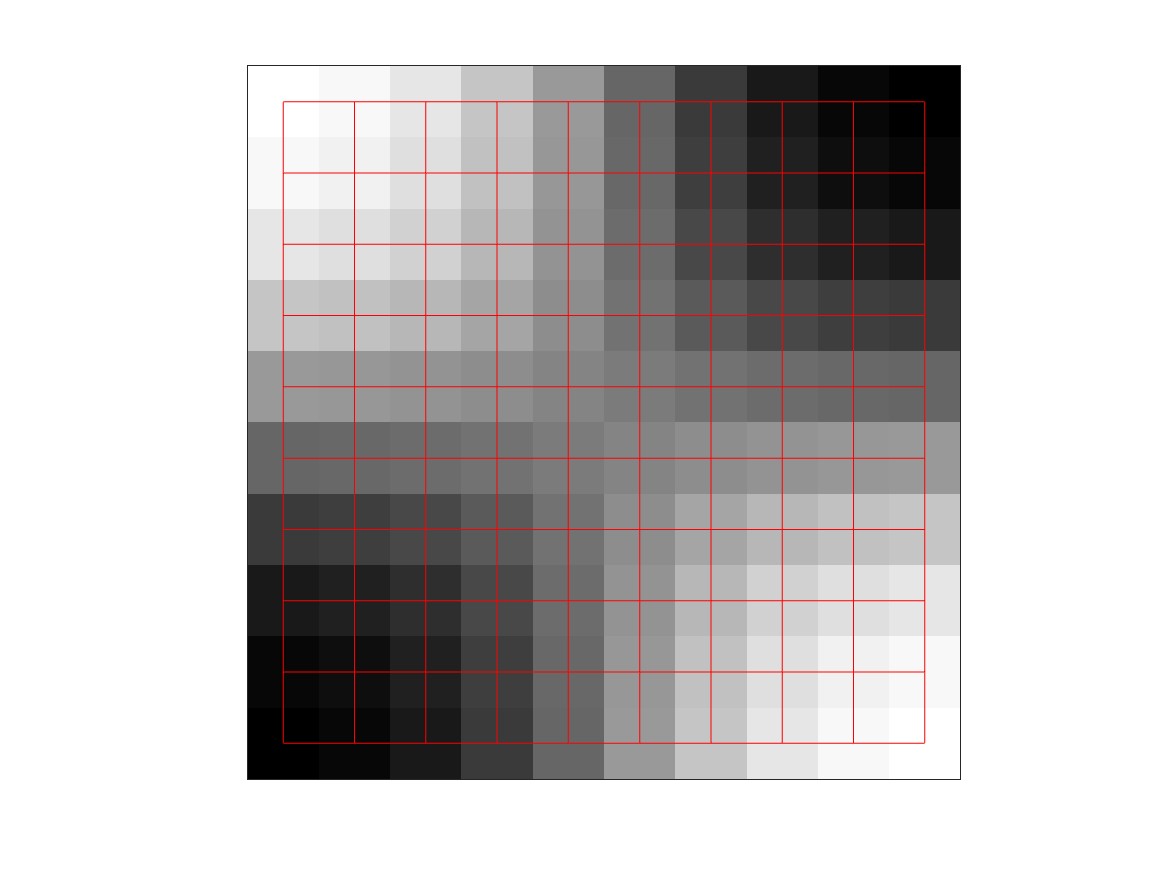} 
  \end{minipage} 
 \vspace{-.5cm}
     \caption{\small Deformation of the images on top  of Fig. \ref{fig2bw}, after the application of 
     the scheme (\ref{schemaFD}) when $T=2$, $\Delta t=0.2$ and  $a=5$.} 
     \vspace{.3cm}
 \label{fig2bwPM}
\end{figure}

 %
 For the treatment of the boundary points, we take into account a discrete version of Neumann boundary conditions. Since the approximation method is explicit
 in time, suitable restrictions on $\Delta t$ should be imposed for stability.
 For the details we address to \cite{Perona}. This last issue is not relevant to the discussion
 of our scheme (actually, the parameter $\gamma$ in  (\ref{step2D})  can be either positive or negative), since the aim is to
 implement the method just for a few iterations.

There are clear differences with respect to the scheme proposed in this paper.
In particular, the second order operator in (\ref{perona}) is straightly approximated on the grid points with integer indices. This realizes
an approximation error only of the second order in the space variable.
As we said in section 3, the anisotropy coefficient of the Laplace operator acts differently in our case,
suggesting  a reverse of the arrow of time.
Finally, the Perona-Malik scheme does not guarantee the preservation of the invariant situations 
examined in 
the previous section. For instance, in  Fig.  \ref{fig2bwPM}, we show what happens by applying 
(\ref{schemaFD}) to the images on top of Fig.  \ref{fig2bw}.

\section{Tests on complex images}\label{sec5}

A typical application in image processing is edge detection, where the 
techniques already available are numerous. 
Remaining in the context of anisotropic diffusion, after applying
our procedure, the resulting image may be treated by
selecting the pixels that tend to diverge from the average.
In practice, the new image is renormalized between  1 and 256 and
the following filter is applied:
\begin{equation}\label{taglio}
{\cal F}_{\tau}(\phi )=\left\{ \begin{array}{l}
256, \quad \phi \ge\tau \ \ {\rm or} \quad \phi \le 256-\tau,\\[3mm] 
1,\ \  {\rm otherwise},\
\end{array} \right. 
\end{equation}
where $\tau$ is a suitable integer between 128 and 256.

%
%
%



\begin{figure}[h!] 
 \vspace{-.4cm}
\hspace{-1.6cm}
  \begin{minipage}[b]{0.59\linewidth}
    \includegraphics[width=1.1\linewidth ]{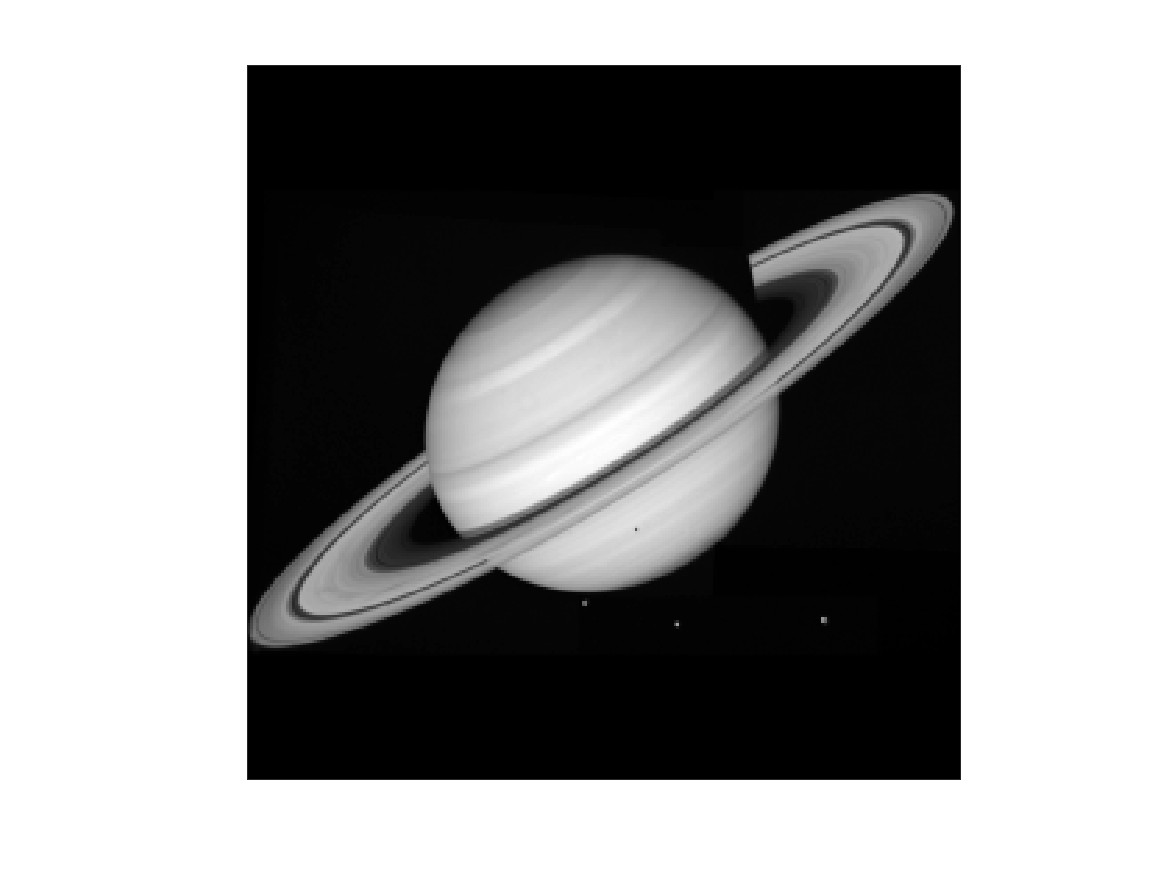} 
  \end{minipage} 
   \hspace{-.8cm}
  \begin{minipage}[b]{0.59\linewidth}
    \includegraphics[width=1.1\linewidth]{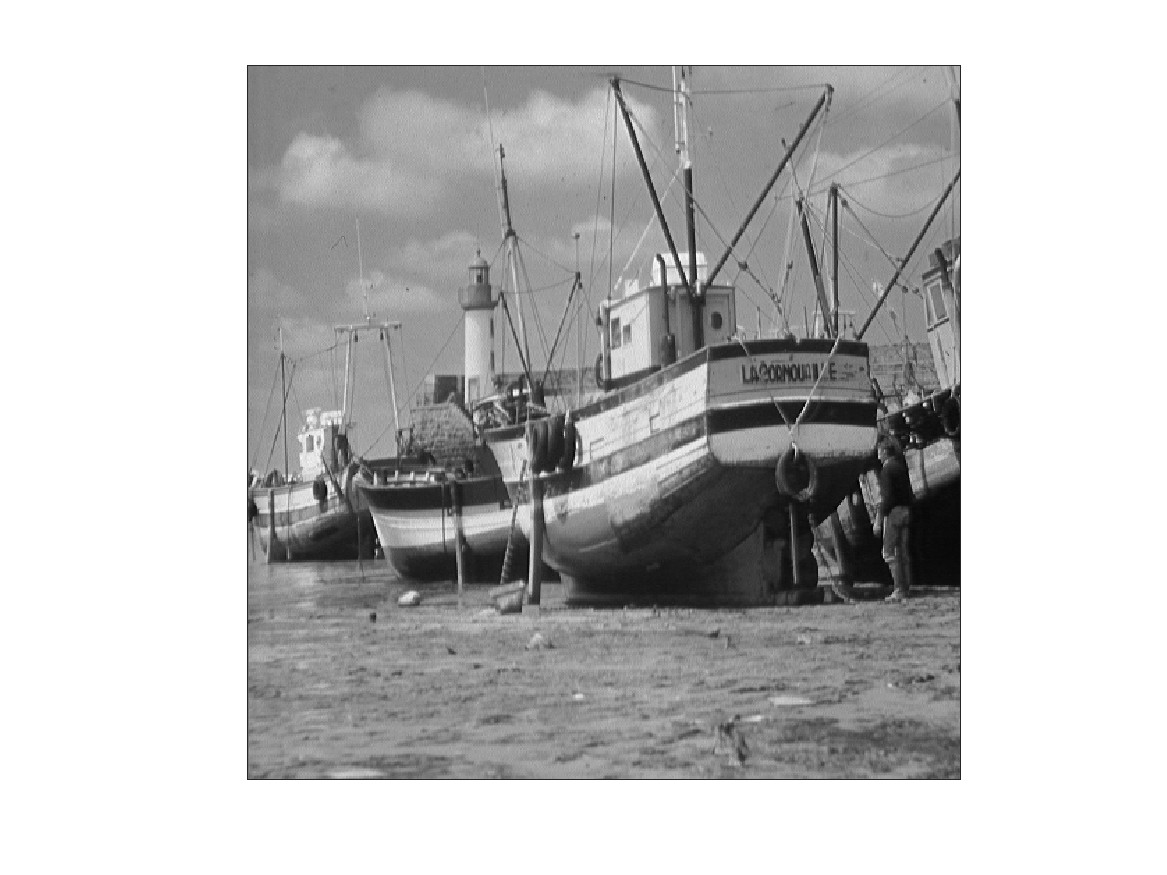} 
  \end{minipage} 
  \vspace{-1.5cm}
     \caption{\small The two images used for the tests.}
 \label{figRealImage} 
\end{figure}


\begin{figure}[ht!] 
\hspace{-1.7cm}
  \begin{minipage}[b]{0.6\linewidth}
    \includegraphics[width=1.1\linewidth ]{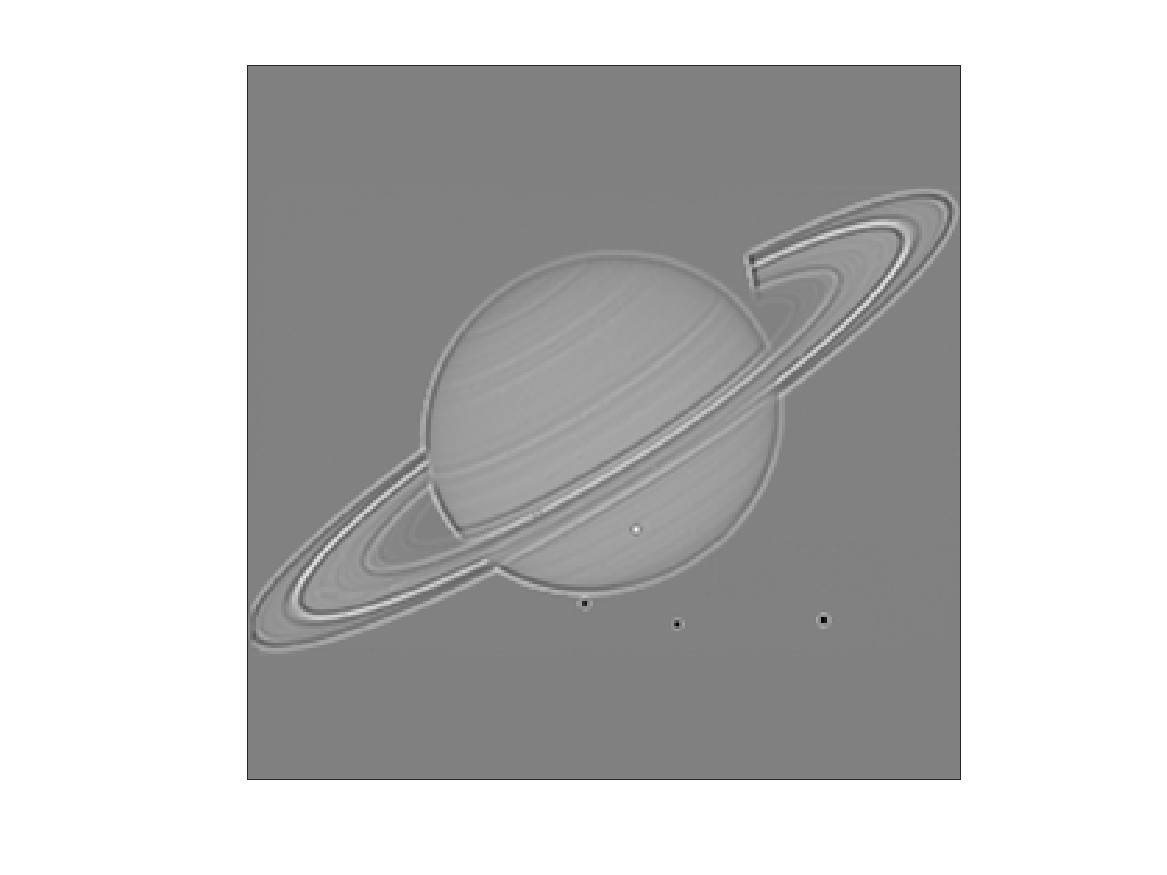} 
  \end{minipage} 
   \hspace{-.8cm}
  \begin{minipage}[b]{0.6\linewidth}
    \includegraphics[width=1.1\linewidth]{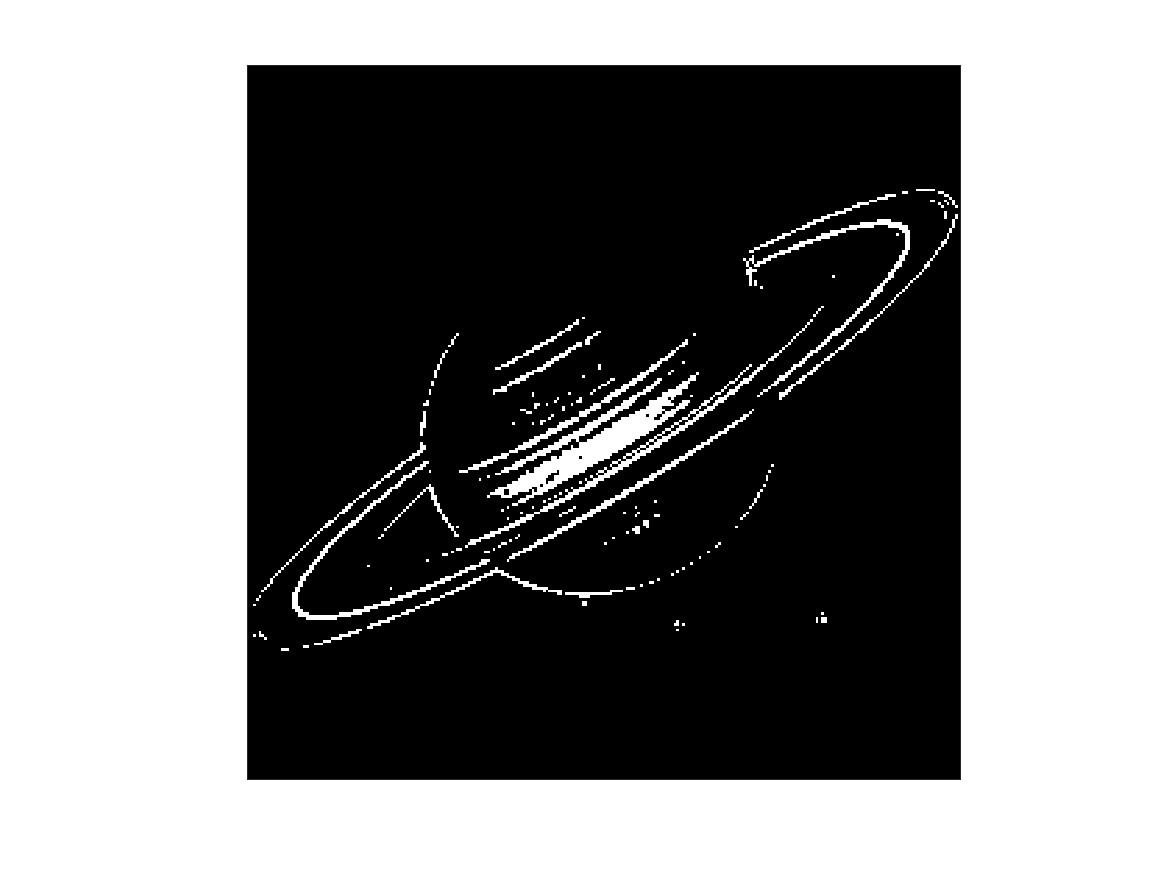} 
  \end{minipage} 
 \vskip-1truecm
  \begin{minipage}[b]{0.6\linewidth}
  \hspace{-1.7cm}
    \includegraphics[width=1.1\linewidth]{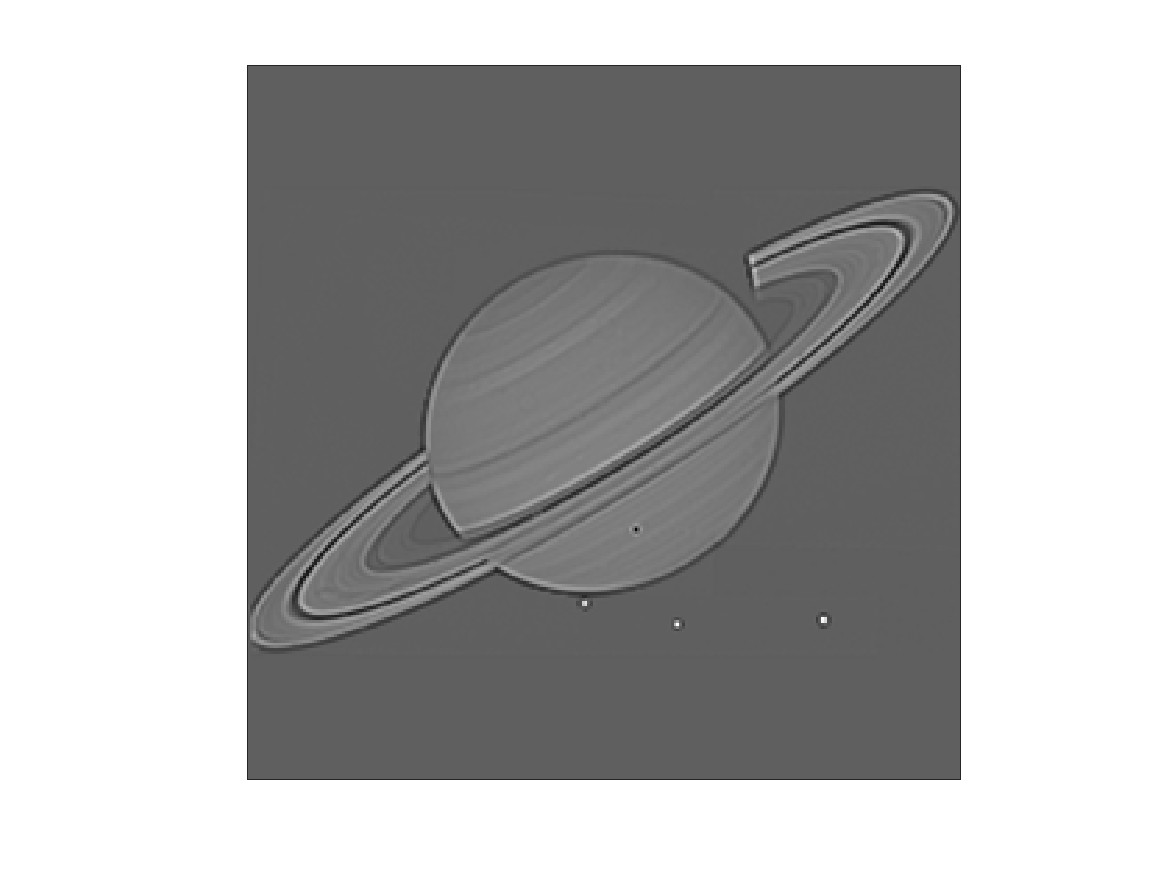} 
  \end{minipage}
 \hspace{-2.4cm}
  \begin{minipage}[b]{0.6\linewidth}
    \includegraphics[width=1.1\linewidth]{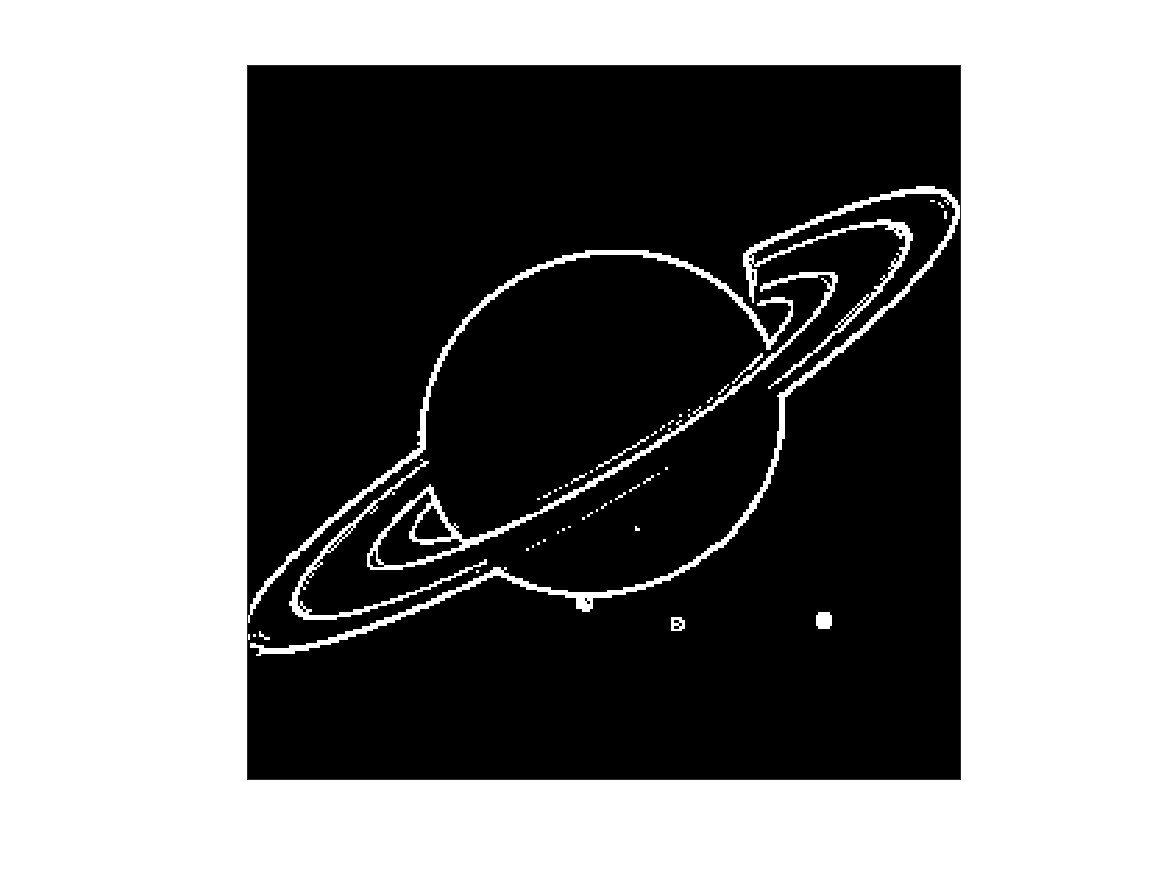} 
  \end{minipage}
  \vspace{-1.cm}
     \caption{  \small The Saturn image has been treated with one step of the algorithm (\ref{step2D})
     ($\gamma=8$ top, $\gamma=-8$ bottom) and a passage of the cut off procedure (\ref{taglio}) with $\tau=162$.
     } 
 \label{figRealImageSaturn} 
\end{figure}

\begin{figure}[ht!] 
 \begin{minipage}[b]{0.6\linewidth}
  \hspace{-1.7cm}
    \includegraphics[width=1.1\linewidth]{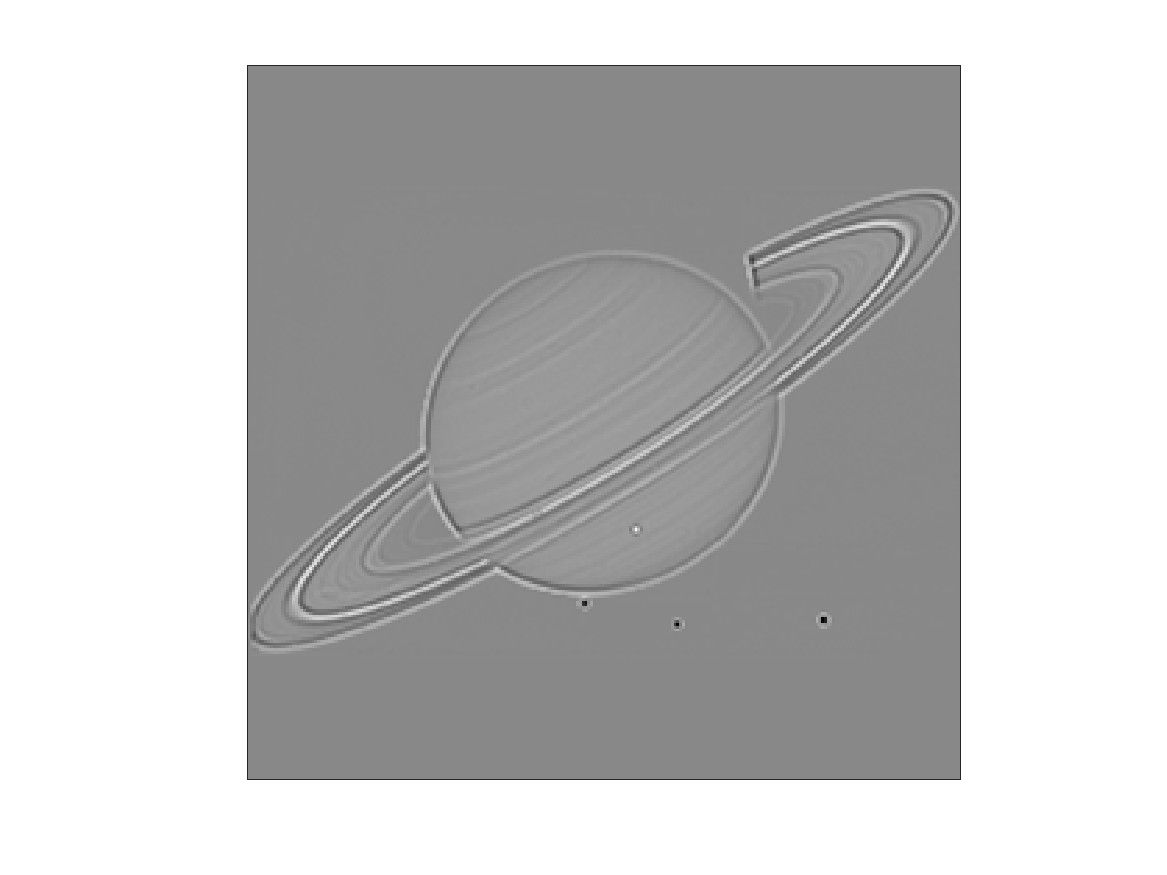} 
  \end{minipage}
 \hspace{-2.4cm}
  \begin{minipage}[b]{0.6\linewidth}
    \includegraphics[width=1.1\linewidth]{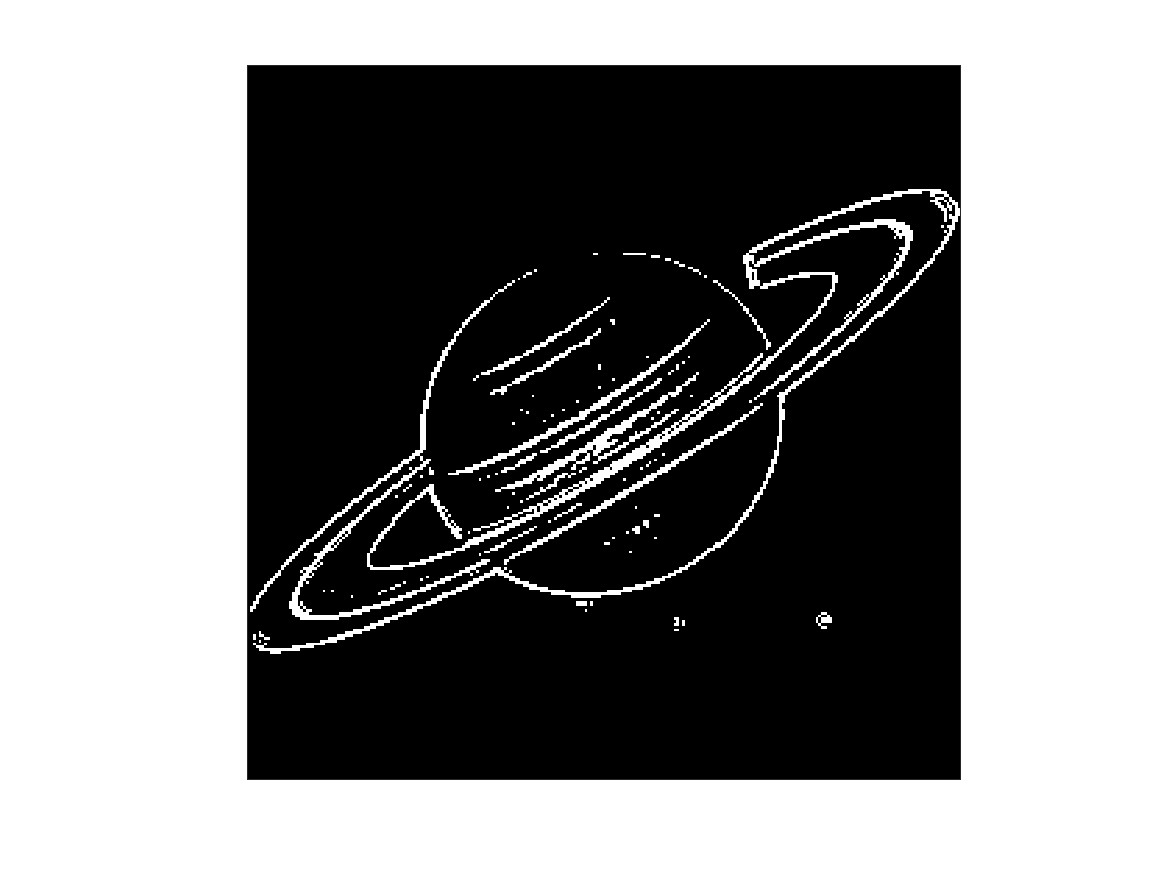} 
  \end{minipage}  
 \vskip-1truecm  
   \begin{minipage}[b]{0.6\linewidth}
  \hspace{-1.7cm}
    \includegraphics[width=1.1\linewidth]{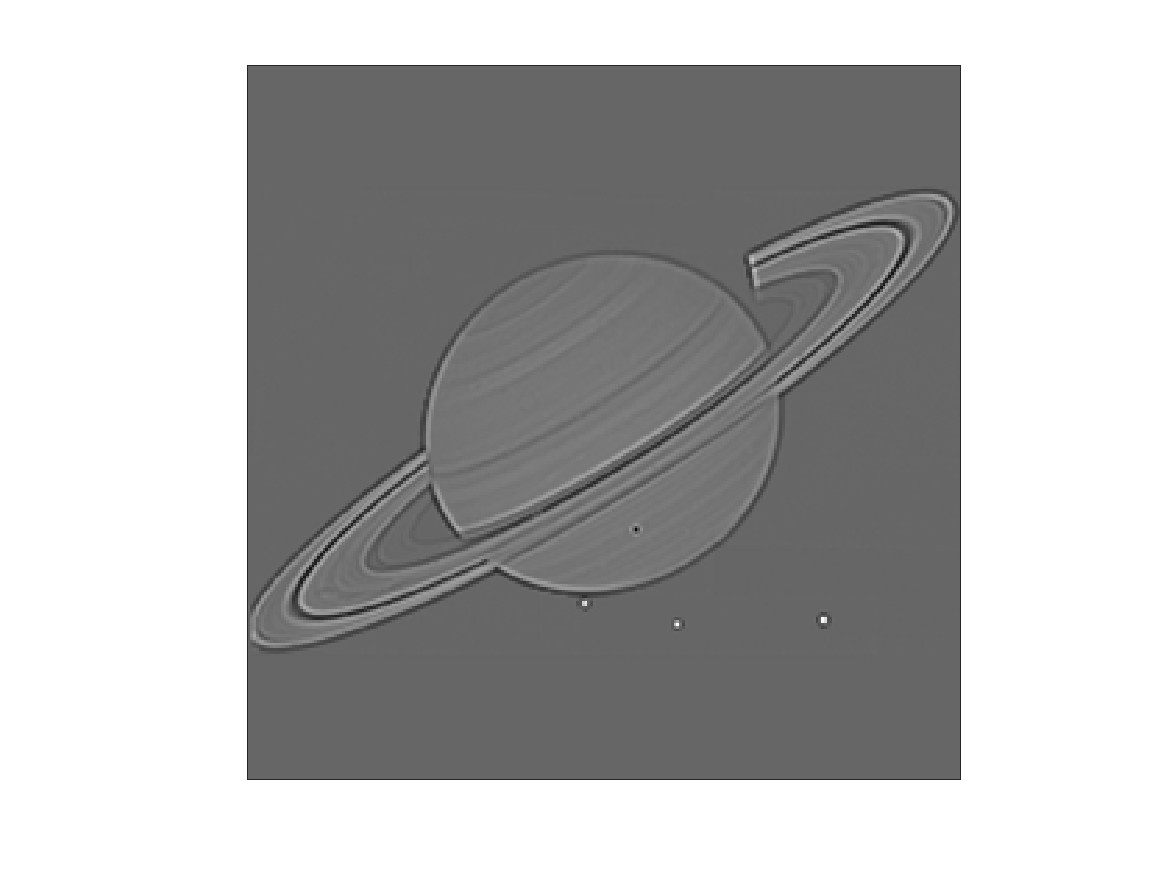} 
  \end{minipage}
 \hspace{-2.4cm}
  \begin{minipage}[b]{0.6\linewidth}
    \includegraphics[width=1.1\linewidth]{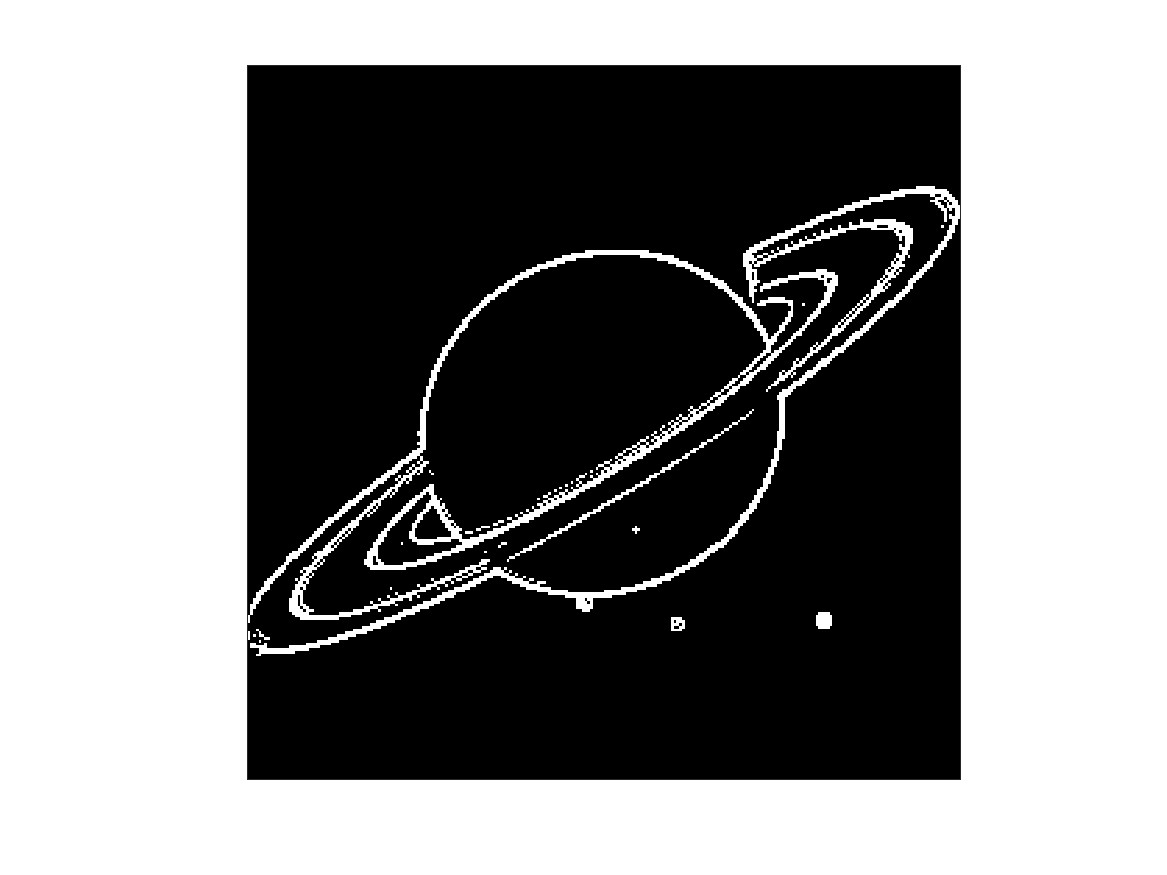} 
  \end{minipage} 
  \vspace{-1.cm}
     \caption{  \small The Saturn image has been treated with one step of the algorithm (\ref{step2D})
     ($\gamma=15$ top, $\gamma=-15$ bottom) and a passage of the cut off procedure (\ref{taglio}) with $\tau=155$.} 
 \label{figRealImageSaturn2} 
\end{figure}

\begin{figure}[ht!] 
\hspace{-1.7cm}
  \begin{minipage}[b]{0.6\linewidth}
    \includegraphics[width=1.1\linewidth]{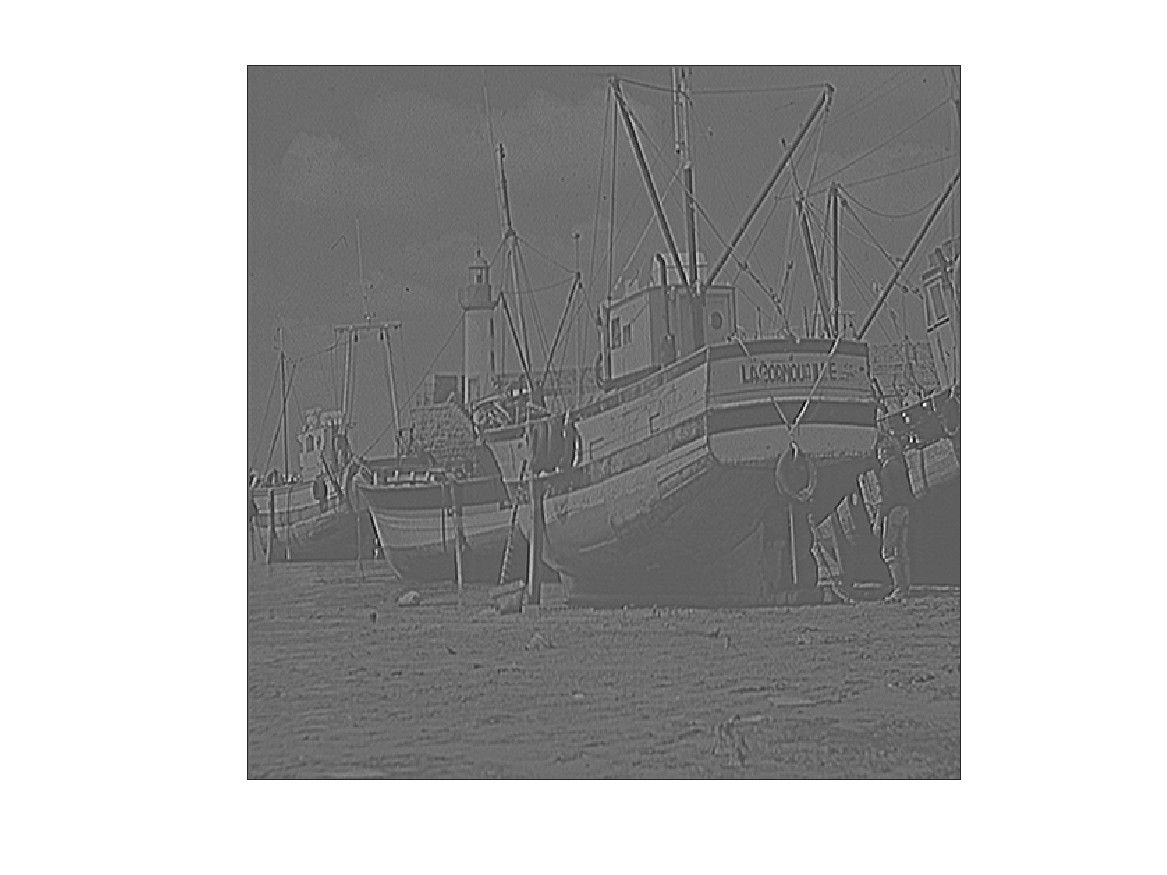} 
  \end{minipage}
 \hspace{-.8cm}
  \begin{minipage}[b]{0.6\linewidth}
    \includegraphics[width=1.1\linewidth]{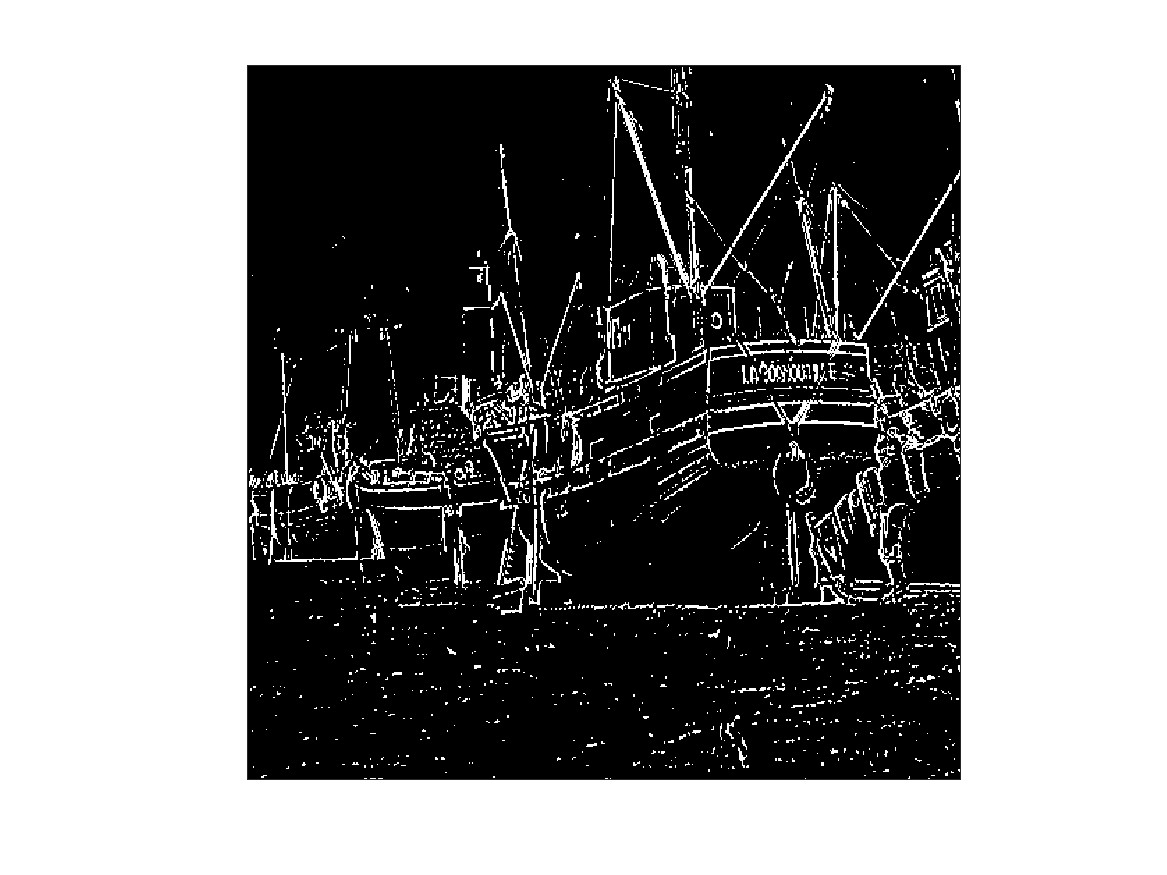} 
  \end{minipage} 
\vspace{-1.cm}
     \caption{  \small The image of the boat has been treated with one step of the algorithm (\ref{step2D})
     ($\gamma=-8$) and a passage of the cut off procedure (\ref{taglio}) with $\tau=165$.
      } 
 \label{figRealImageShip} 
\end{figure}

We  carried out some tests on  the two images shown in  Fig. \ref{figRealImage}. 
We mainly worked on edge detection, although this is only one of the possible applications.
After running a single step of the  algorithm, we get the outcomes displayed 
in  Figs. \ref{figRealImageSaturn}, \ref{figRealImageSaturn2}, \ref{figRealImageShip}.
Different values of $\gamma$ (positive and negative) in the scheme (\ref{step2D}), and $\tau$ in the filtering procedure (\ref{taglio}), are chosen for
comparisons.
 The results are not bad at all, especially
if we consider that they have been obtained after just one step of the iterative procedure (\ref{step2D}),
thus with a cost proportional to the number of pixels. 
We point out that similar algorithms (such as the Perona-Malik, described in the previous
section) require instead a large number of iterations. 
As expected, the use of negative values of $\gamma$ (i.e., by going backward in time) provides
sharper details. Of course, the goodness of the performances should be judged depending on
the type of target one has in mind.
We also tried some experiments (not shown here) by iterating (\ref{step2D}) more times.
The results are not so neat, however, 
since the images start becoming too blurred.

\begin{figure}[ht!] 
\hspace{-1.7cm}
  \begin{minipage}[b]{0.6\linewidth}
    \includegraphics[width=1.1\linewidth ]{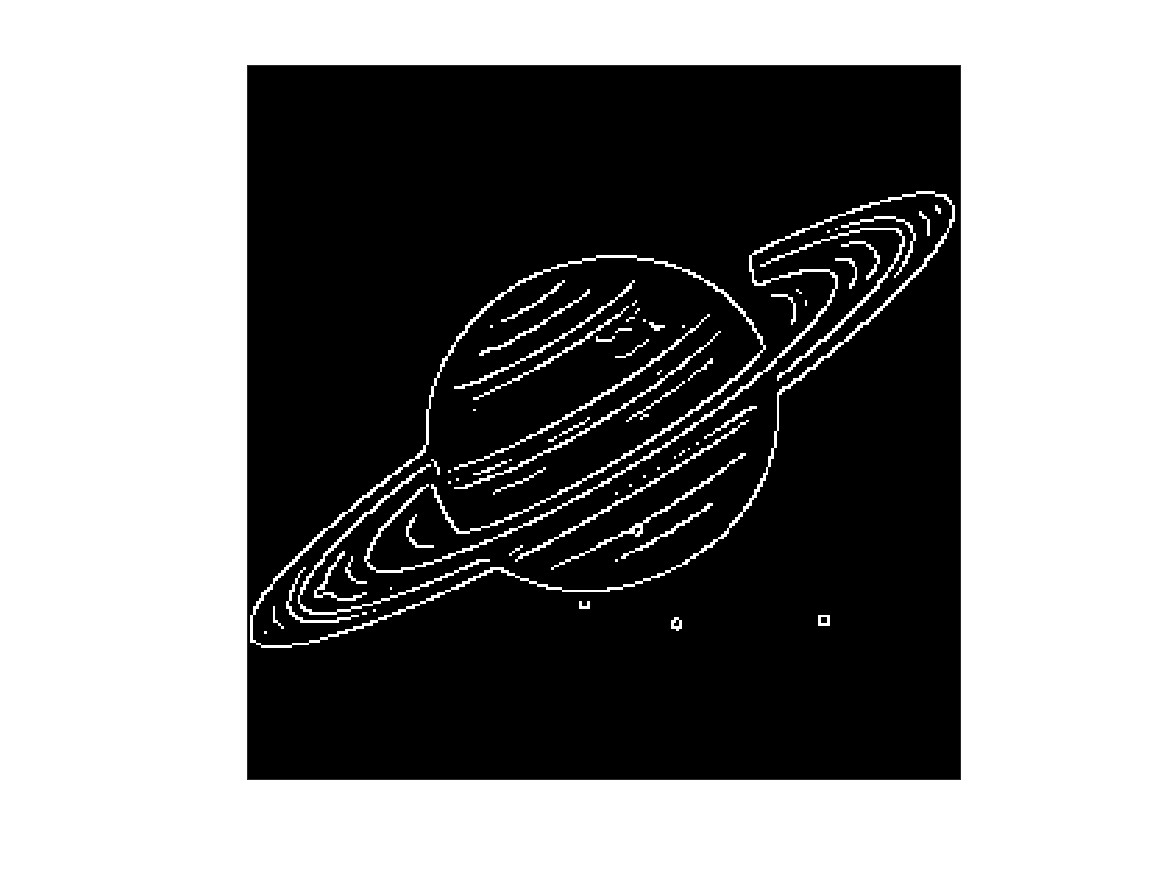} 
  \end{minipage} 
   \hspace{-.8cm}
  \begin{minipage}[b]{0.6\linewidth}
    \includegraphics[width=1.1\linewidth]{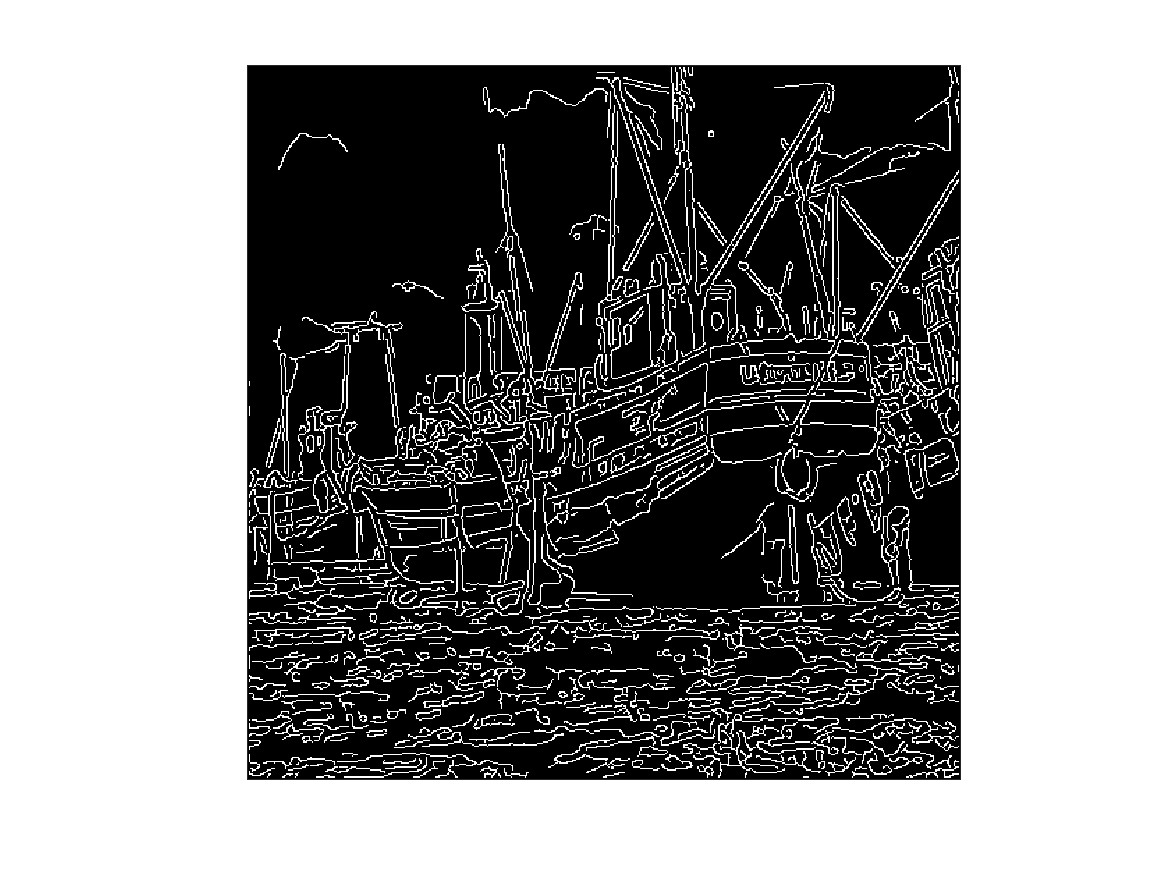} 
  \end{minipage} 
  \vspace{-1.cm}
     \caption{\small The results of the Canny algorithm applied to the two images used for the tests.}
 \label{figCannyImage} 
\end{figure}

Finally, we conclude this section by showing in Fig. \ref{figCannyImage}  the results obtainable 
with a standard algorithm, such as the Canny edge detector \cite{Canny}. In particular the images in  Fig. \ref{figCannyImage} are obtained by using  the Matlab command {\tt edge}. By default,  {\tt edge} uses the Sobel  detection method. We refer to the MathWorks documentation for further insights.

\section{Conclusions}\label{conclusions}

A suitable diffusive PDE with anisotropic coefficient has been proposed, together with 
a high-order space discretization making use of staggered grid points. Some exploratory tests
have been carried out in the framework of  digital image processing. A single
step of the algorithm is already sufficient to emphasize the contour lines of the objects.
Thus, the procedure can successfully find application in edge detection or denoising techniques, also
in conjunction with other well-known methods.
\section*{Acknowledgments}
Both authors are members of GNCS-INDAM (Gruppo Nazionale Calcolo Scientifico - Istituto Nazionale di Alta Matematica), Rome. 

\end{document}